\documentclass[11pt]{article}

\usepackage[latin1]{inputenc}
\usepackage{amsmath,amssymb, bbm}
\usepackage{latexsym}

\newtheorem{theorem}{Theorem}[section]
\newtheorem{corollary}[theorem]{Corollary}
\newtheorem{lemma}[theorem]{Lemma}
\newtheorem{proposition}[theorem]{Proposition}
\newtheorem{definition}[theorem]{Definition}

\newtheorem{remark}[theorem]{Remark}
\numberwithin{equation}{section}

\def\dist{\operatorname{dist}}

\def\il{\int\limits}

\def\vp{{\varphi}}
\def\qq{{\qquad}}

\def\wt{\widetilde}
\def\R{{\mathbb R}}

\def\E{{\mathbb E}}
\def\P{{\mathbb P}}
\def\H{{\mathbb H}}
\def\Q{{\mathbb Q}}

\def\eps{\varepsilon}

\def\al{\alpha}

\def\phi{\varphi}

\def\lam{{\lambda}}
\def\ol{\overline}

\def\proof{{\medskip\noindent {\bf Proof. }}}

\def\qed{{\hfill $\square$ \bigskip}}

\def\dist{{\mathop {{\rm dist\, }}}}

\def\square{{\vcenter{\vbox{\hrule height.3pt
        \hbox{\vrule width.3pt height5pt \kern5pt
           \vrule width.3pt}
        \hrule height.3pt}}}}

  \def\sC {{\cal C}}
 \def\sE {{\cal E}} \def\sF {{\cal F}}
  
\def\sJ {{\cal J}}  \def\sL {{\cal L}}
 \def\sN {{\cal N}}

 \def\bE {{\mathbb E}}

\def\bP {{\mathbb P}}  \def\bR {{\mathbb R}}

\def\nn{\nonumber}

\def\q{\quad}
\def\half{{\textstyle \frac12}}

\def\ignore#1{}

\def\pf{\noindent{\bf Proof.} }

\def\a4xi{{(A4)($\xi$)}}

\begin{document}

\title{\bf Non-local Dirichlet Forms and Symmetric Jump Processes}

\author{
Martin T. Barlow\footnote{Research partially supported by NSERC
(Canada)}, \ Richard F. Bass\footnote{Research partially supported
by NSF grant DMS-0601783.}, \
Zhen-Qing Chen\footnote{Research partially supported by NSF grant
DMS-06000206.}, \ and Moritz Kassmann\footnote{Research partially
supported by DFG (Germany) through Sonderforschungsbereich 611.} }


\maketitle


\begin{abstract}
\noindent We consider the symmetric
non-local Dirichlet form $(\sE, \sF )$ given by
\[ \sE (f,f)=\il_{\R^d} \il_{\R^d} (f(y)-f(x))^2 J(x,y) \, dx\, dy  \]
with $\sF$ the closure of the set of $C^1$ functions on $\R^d$ with
compact support with respect to $\sE_1$,
where $\sE_1 (f, f):=\sE (f, f)+\int_{\R^d} f(x)^2 dx$, and
where the jump kernel $J$ satisfies
\[ \kappa_1|y-x|^{-d-\alpha} \leq J(x,y) \leq \kappa_2|y-x|^{-d-\beta} \]
 for $0<\alpha< \beta <2, \, |x-y|<1$.
This assumption allows the corresponding jump process to
have jump intensities whose size depends on the position of the process
and the direction of the jump.  We prove upper and lower
estimates on the heat kernel.
We construct a strong Markov process corresponding to $(\sE, \sF)$.
We prove a parabolic Harnack inequality for nonnegative
functions that solve the heat equation with respect to $\sE$.
Finally we construct an example where the corresponding harmonic functions
need not be continuous.


\end{abstract}

\section{Introduction}

In this paper we introduce a class of symmetric Markov processes of
pure jump type.
Our assumptions allow the jump intensities to depend on both the
position of the process and the direction of the jump. Thus our
processes can be highly anisotropic.
Although very little regularity is assumed, nevertheless
we are able to obtain a number of results concerning these processes.

We begin by considering symmetric non-local Dirichlet forms.
Set
\begin{eqnarray}\label{Dirform}
 \sE (f,f)&=&\il_{\R^d} \il_{\R^d} (f(y)-f(x))^2 J(x,y) \, dx\, dy \;, \\
\label{eqn:domain}
\sF&=& \overline{C^1_c (\R^d)}^{\sE_1}  \; ,
\end{eqnarray}
where the jump kernel
$J(x,y)$ is a function of $x$ and $y$ satisfying the following conditions
\begin{itemize}
\item[(A1)] $J(x,y) = J(y,x)$  for all $x$ and $y$;
\item[(A2)] $J(x,y) = 0$  for $|x-y| \geq 1$;
\item[(A3)] There exist $\alpha, \beta \in (0,2)$, $\beta > \alpha $
and positive $\kappa_1, \kappa_2$ such that
\[ \kappa_1 |y-x|^{-d-\alpha} \leq J(x,y) \leq \kappa_2 |y-x|^{-d-\beta}
\qquad \hbox{for } |y-x|<1 \;. \]
\end{itemize}
Here   $\sE_1(f, f):= \sE(f,f)  +\|f\|_2^2$,
  $C^1_c(\R^d)$ denotes the space of $C^1$ functions on $\R^d$ with
compact support, and $\sF$ is the closure of $C^1_c(\R^d)$ with respect
to the metric $\sE_1(f,f)^{1/2}$.

\medskip

We obtain the following four main results in this paper.  We emphasize
that we make no continuity assumptions whatsoever on the jump kernel
$J$.

\begin{description}

\item{(i)}
We show that there is a Hunt $X$ process associated to $(\sE, \sF)$, and
$X$ has a symmetric transition density function $p(t,x, y)$, with
respect to Lebesgue measure on $\R^d$.
 We
call this transition density function the heat kernel corresponding
to $(\sE, \sF)$ and derive upper and lower bounds.

\item{(ii)} We show that the strong Markov process $X$ corresponding to $(\sE, \sF)$
is conservative and can
be realized as a weak limit of certain more regular jump processes.

\item{(iii)} We establish a parabolic Harnack inequality for nonnegative functions
that solve the heat equation with respect to $(\sE, \sF)$.

\item{(iv)}  We construct a counterexample to show that harmonic functions with
respect to $X$ need not be continuous on $\R^d$.

\end{description}

\noindent We now discuss each of these points in turn.

\medskip
We first show there exists a Hunt process associated to $(\sE, \sF)$.

\begin{theorem}\label{thmzero}
Suppose (A1)--(A3) hold.
There exists  $\sN \subset \R^d$ having zero
capacity with respect to the Dirichlet
form $(\sE, \sF)$ and there is a Hunt process $(X, \P^x)$ with state space
$\R^d \setminus \sN$ such that for every $f\in L^2(\R^d, dx)$ and $t>0$,
$x\mapsto \E^x [ f(X_t)]$ is a quasi-continuous version of $T_t f$,
where $\{T_t, t\geq 0\}$ is the $L^2$-semigroup associated with the closed form $(\sE, \sF)$.
\end{theorem}

\bigskip

 Note that the Hunt process $X$ can start from any point in
$\R^d \setminus \sN$ and that once it starts from $\R^d\setminus
\sN$ the process $X$ together with its left hand limits takes values
in $\R^d \setminus \sN$ up to and strictly before  its lifetime
$\zeta$. The set $\sN$ is called the properly exceptional set of $X$
(or, equivalently, of $(\sE, \sF)$) and it has zero Lebesgue
measure. For simplicity, sometimes we just say that $X$ is a Hunt
process associated with $(\sE, \sF)$ starting from quasi-everywhere
in $\R^d$. For more on terminology and properties of Dirichlet
forms, we refer the reader to \cite{FOT94}.

Let $P(t,x,dy)$ be the transition probability for the Hunt process $X$
associated with $(\sE, \sF)$.

\begin{theorem}\label{rough-above-lem}
Assume conditions (A1)--(A3) hold. There is a properly exceptional
set $\sN \subset \R^d$ of $X$, a positive symmetric kernel
$p(t,x,y)$ defined on $(0, \infty)\times (\R^d \setminus \sN) \times
(\R^d \setminus \sN)$, and positive constants $C_1$ and $C_2$
(depending on the constants in (A1)--(A3))
 such that $P(t,x,dy) = p(t,x,y) dy$, and
\begin{align}
p(t, x, y) \leq C_1 t^{-d/\alpha} e^{C_2t}
\qquad \hbox{for every } t>0 \hbox{ and } x, y \in \R^d\setminus \sN.
\label{rough-above}
\end{align}
Moreover, for every $t>0$ and $y\in \R^d\setminus \sN$,
$x\mapsto p(t, x, y)$ is quasi-continuous on $\R^d$.
\end{theorem}

\medskip
We also obtain lower bounds on the heat kernel.  Let $B(x,r)$ denote
the open ball of radius $r$ centered at $x$, and $p^B(t,x,y)$ be the
transition densities for the subprocess of $X$ killed upon exiting
the ball $B$.

\begin{theorem}\label{lowerbound-theo}
Assume conditions (A1)--(A3) hold. Let  $y_0\in \R^d$, $T>1/2$, and
$\delta\in (0,1/2)$. Let $R>0$ and $B=B(y_0,R)$. There exists a
properly exceptional set $\sN$ and a positive constant $C$ that
depends on $R, T, \al, \beta, \kappa_1, \kappa_2$, and $\delta$, but
not on  $y_0$ such that    for all $t\in [\delta, T]$
\begin{equation} \label{lowerbound}
p^B(t,x,y) \ge C  \end{equation}
for every $(x,y)\in (B(y_0, 3R/4)\setminus \sN)\times
 (B(y_0, 3R/4)\setminus \sN)$.
\end{theorem}

\begin{remark}\label{R:1.4}
{\rm The jump kernel $J$ does not have any scaling properties, and
so one should not expect the results in Theorems \ref{rough-above-lem} and
\ref{lowerbound-theo} to be scale invariant. In particular, the
constant $C$ in Theorem \ref{lowerbound-theo} depends on $R$. }
\end{remark}

\medskip
One of the difficulties in working with the process associated
to $(\sE, \sF)$ is that we do not know that
$\sF$ defined above is equal to
\[ \sF_\infty=\{f\in L^2(\R^d,dx):  \sE(f,f)<\infty\}. \]
(A similar problem arises when studying the minimal Brownian motion on
a manifold.)  To circumvent this and other difficulties, we
will approximate $J$ by jump kernels $J_\xi$ which have regular
behavior for $|x-y|\le \xi$.
We introduce the following condition, which we will assume from time to time.
Let $\xi >0$.

\medskip \noindent \a4xi \quad The jump kernel $J$ satisfies
\begin{align*}
 J(x,y)=\kappa_2|y-x|^{-d-\beta} \qquad \hbox{when } \, |y-x|<\xi.
\end{align*}
If $J$ satisfies (A1)--(A3), define
\begin{equation}\label{Jxidefine}
 J_\xi(x,y)=J(x,y)\mathbbm{1}_{(|x-y|>\xi)}
 + \kappa_2|x-y|^{-d-\beta}\mathbbm{1}_{(|x-y|\leq \xi)}.
\end{equation}
Let $(\sE^{(\xi )}, \sF^{(\xi)} )$ be the regular Dirichlet form on $\R^d$ defined
by (\ref{Dirform})--(\ref{eqn:domain}) with $J_\xi$ in place of $J$.
We can then prove (see Lemma \ref{init-m-bounds} below) that
$$ \sF^{(\xi)}=\{f\in L^2(\R^d, dx): \ \sE^{(\xi )}(f,f)<\infty\}. $$
We also have

\begin{theorem}\label{moscoconv}
As $\xi \to 0+$,
$(\sE^{(\xi )}, \sF^{(\xi)})$ converges to $(\sE, \sF)$ in the sense of Mosco.
\end{theorem}

\begin{remark}\label{R:1.6}{\rm
See Definition \ref{D:2.1} for the definition of Mosco
convergence.

\begin{description}

\item{(i)} Mosco convergence (see \cite{Mosco}) implies
that the semigroups of the processes $X^{(\xi)}$ associated to
$\sE^{(\xi )}$ converge in $L^2(\R^d, dx)$ to the semigroup of the
process $X$. We establish  Theorem \ref{lowerbound-theo} by first
proving the result for $X^{(\xi)}$, with constants independent of
$\xi$; taking a limit then gives the result for $X$.

\item{(ii)} Studying $X$ by first assuming \a4xi and then taking limits
is analogous to a common procedure in the study of elliptic
operators in divergence form. There one often first assumes the coefficients
are smooth and obtains estimates that do not depend on the smoothness,
and then uses a limiting procedure.

\item{(iii)} We prove Theorem \ref{moscoconv} by first establishing a
simple sufficient criterion for Mosco convergence to hold. This is of
independent interest.

\item{(iv)}
It seems to be difficult to establish a similar result if we approximate $J(x, y)$
from below when $|x-y|$ is small.

\item{(v)} As remarked above, we do not know in general that $\sF=\sF_\infty$.
However if the jump intensity kernel $J$ is ``translation
equivalent" near the diagonal, that is, if there exist constants
$c>1$ and $\delta>0$ such that
$$ c^{-1} J(x, y) \leq J(x-z, y-z) \leq c J(x, y) $$
for a.e.  $x, y, z \in \R^d$   with $0<|x-y|< \delta$, then it is
not difficult to show that $\sF=\sF_\infty$.
Since such a result will not be used in this paper, we omit its
proof.
\end{description} }\end{remark}

\medskip

Let $X$ be the Hunt process associated with the regular
Dirichlet form $(\sE, \sF)$ on
$\R^d$, which   has an exceptional set $\sN$.
 We say a function
$h:\R^d\to \R$ is {\em harmonic} on a ball $B(y,r)$ if
$h(X_{t\land \tau_{B(y,r_1)}})$ is a $\P^x$-martingale with right continuous
paths for every $r_1\in (0, r)$ and
every $x\in B(y, r_1)\setminus \sN$. Here
 $\tau_{B(y,r_1)}=\inf\{t: X_t\notin B(y,r_1)\}$.

Set $V_t=V_0-t$ and let $\bP^{(s,x)}$ be the law of $(V_t, X_t)$
started at $(s,x)$.  We say a function $u: [0,\infty)\times \R^d\to \R$
is {\em caloric} on $Q=(a,b)\times B(x_0,r)$ with respect to $\sE$ if
$u(V_{t\land \tau_Q}, X_{t\land \tau_{Q_1}})$ is a
$\bP^{(s,x)}$-martingale with right continuous paths for every open
subset $Q_1$ of $Q$ with $\overline Q_1 \subset Q$ and for every $(s,x)\in
Q_1 \cap \left( \R_+\times (\R^d\setminus \sN) \right)$.  Here
$\tau_{Q_1}=\inf\{t: (V_t,X_t)\notin Q_1\}$.

We prove that nonnegative functions that are caloric with respect to
$\sE$ satisfy a parabolic Harnack inequality.

\begin{theorem} \label{Maintheorem}
Suppose the Dirichlet form $(\sE, \sF)$ is given by
(\ref{Dirform})--(\ref{eqn:domain}) with  $J(x,y)$ satisfying
(A1)--(A3). Let $t_0\geq 0, R\geq 1$, and $T>0$.There exists a
positive real $C=C(\alpha, \beta, \kappa_1, \kappa_2, d, R, T)$ such
that if $x_0\in \R^d$ and $u$ is nonnegative and bounded in $(t_0,
t_0+5T) \times \R^d$ and is caloric on $Q=(t_0, t_0+5T)\times
B(x_0,4R)$ with respect to $\sE$, then
\[ \mbox{\rm ess sup}_{Q^-}\, u \leq C\,\mbox{\rm ess inf}_{Q^+}\, u, \]
where $Q^-=[t_0+T,t_0+2T]\times B(x_0,R)$ and $Q^+=[t_0+3T,t_0+4T]\times B(x_0,R)$.
\end{theorem}

\bigskip

\begin{remark}\label{R:1.8} {\rm Concerning the hypotheses and statement of
Theorem \ref{Maintheorem}, we make the following remarks.

\begin{description}

\item{(i)} We assume $u(t, \cdot)$ is bounded in $\R^d$ in order to
ensure that the random variable $u(V_{t\land \tau_Q}, X_{t\land
\tau_{Q_1}})$ is integrable.
 However the constant $C$ does not depend on this bound.

\item{(ii)} Harmonic functions are caloric, so the parabolic Harnack
inequality implies that an elliptic Harnack inequality also holds.

\item{(iii)} Assumption (A3) does not satisfy any type of scaling property.
As a result,  one cannot expect
the parabolic Harnack inequality to be scale invariant, i.e., that the
constant $C$ can be chosen independently of $R$ or $T$. Since an
example in \cite{BaKa05a} shows that scale invariance can fail for
the elliptic Harnack inequality, it can also fail for the parabolic
Harnack inequality. This phenomenon is well-known in the theory of degenerate partial differential equations, see \cite{ChWh86}, \cite{GuWh90}.

\item{(iv)} We shall see in Theorem \ref{Theo-counterexample} below
that it may not be possible to extend a harmonic function $h$ to a
continuous function on $\bR^d$.  Thus we have to use the essential
supremum and essential infimum in Theorem \ref{Maintheorem}.

\item{(v)} Assumption (A2) rules out
any jumps of size larger than 1. An example in \cite{BaKa05a} shows
that the large jumps, although in many ways less interesting, can cause the
Harnack inequality to fail.
For similar reasons we cannot replace the ball of radius 1 by
arbitrarily small balls in Theorem \ref{Maintheorem}.

\item{(vi)} We allow $0<\alpha<\beta<2$ with no other restriction on $\alpha$
and $\beta$. This should be contrasted with the situation in
\cite{BaKa05a}, which considers non-local operators that are
non-symmetric, and where in addition it was required that
$\beta-\alpha<1$.
\end{description} } \end{remark}

Following \cite{FaSt86} many papers have used heat kernel estimates to
prove Harnack inequalities. The usual procedure is to obtain an
oscillation inequality, and from this
one obtains
a Harnack inequality. We cannot use this approach here, since
our counterexample shows that the constant in the oscillation
inequality can blow up as the radius $r_n$ of the ball
approaches 0.
Instead, the proof of Theorem \ref{Maintheorem} uses a
balayage argument; this approach is new and is of independent interest.

\medskip

As we mentioned above, harmonic functions need not be continuous.
We prove

\begin{theorem}\label{Theo-counterexample} There exists a Dirichlet form
$(\sE, \sF)$ given by (\ref{Dirform})--(\ref{eqn:domain}) with
the jump kernel $J$ satisfying (A1)--(A3), but
where there also exists a bounded harmonic function
that cannot be extended to be a continuous function on $B(0,1)$.
\end{theorem}

\begin{remark}\label{R:1.11}
\begin{description} {\rm
\item{(i)}  We will also show that
continuity can fail for $P_tf$, even when $f$ is smooth.

\item{(ii)} Our construction also gives an example of a martingale problem
for which uniqueness fails.
The example also shows that the existence of  a properly exceptional
set $\sN$ in Theorem \ref{thmzero} and Theorem \ref{rough-above-lem}
is essential and cannot be dropped in general.

\item{(iii)}  The harmonic function we construct may be continuous
outside a set $\sN$ of capacity 0.
}
\end{description}
\end{remark}

Heat kernel estimates and Harnack inequalities have a long history in
the theory of partial differential equations. After path breaking work
by DeGiorgi \cite{DeG57} and Nash \cite{Nas58} on regularity, Moser
\cite{Mos61} proved a scale invariant Harnack inequality for functions
that are harmonic with respect to second order elliptic operators in
divergence form.  This was extended in Moser \cite{Mos64} to solutions
to the heat equation, i.e., the parabolic case; see also \cite{Mos71}.
A quite different proof of this was given in Fabes-Stroock
\cite{FaSt86}.  The Harnack inequality for operators in
nondivergence form was established by Krylov-Safonov \cite{KrSa80}.
However the corresponding theory of Harnack inequalities for
jump processes is still largely unknown.

Non-local operators such as those considered in this paper
arise in the study of models of financial markets (see \cite{SS06}
and the references therein). They
also arise in the study of the Dirichlet-to-Neumann map, particularly
for subelliptic operators or in rough domains.

Harnack inequalities for non-local operators have been considered in
\cite{BaLe02a} and \cite{SoVo04} for fixed order, non-symmetric
operators, \cite{BaLe02b} and \cite{CK03} for fixed order, symmetric
operators. A scale dependent Harnack inequality has been established
in \cite{BaKa05a} for variable order, non-symmetric operators.
Additionally, regularity of harmonic functions is considered for
variable order, non-symmetric operators in \cite{BaKa05b},
\cite{HuKa05}.  For heat kernel estimates and parabolic Harnack
principle for symmetric non-local Dirichlet forms on $d$-sets, see
\cite{CK03}, \cite{HuKu05} for fixed order and \cite{CK06} for
variable order.
See \cite{SchUe} for related results for processes given in terms
of pseudo-differential operators.

\medskip

The paper is organized as follows. In the next section we obtain some
upper bounds for the fundamental solution of the operator
corresponding to $\sE$ and in Section 3 we consider lower
bounds. The Mosco convergence   is proved in Section 4.  The
parabolic Harnack inequality is established in Section 5.
The counterexample is constructed in Section 6.
We use $c_i, c$ or $C$ to denote finite positive constants that depend only on
  $\al$, $\beta$, $\kappa_i$ or $d$ and whose exact value is not
  important and may change from line to line.
 Further dependencies are mentioned explicitly.
We denote the Lebesgue measure of a Borel set $A$ by $|A|$. If $A$ is a
Borel set and $Y$ a right continuous process, we use the notation
\begin{equation}\label{hittingtimedef}
T^Y_A=T_A=\inf\{t>0: Y_t\in A\}, \qquad \tau^Y_A=\tau_A=\inf\{t>0: Y_t\notin A\}.
\end{equation}
For processes $Y$ with paths that are right continuous with left
limits, we let $Y_{t-}$ be the left hand limit at time $t$ and
$\Delta Y_t :=Y_t-Y_{t-}$ the jump at time $t$.

\section{Upper bounds for the heat kernel}

Throughout this section we will assume that the jump intensity
kernel $J$ satisfies (A1)--(A3). We begin with the proof of  Theorem
\ref{thmzero}, which  is easy.
\medskip

\noindent{\bf Proof of Theorem \ref{thmzero}:}
Let $C_\infty (\R^d)$ denote the space
of continuous functions on $\R^d$ that vanish at infinity and let
$\| \cdot \|_\infty$ denote the supremum norm in $C_\infty (\R^d)$.
It is easy to check by using  Fatou's lemma that
the bilinear form $(\sE, \sF)$ is a closed form (cf. \cite[Example 1.2.4]{FOT94}).
As $C^1_c(\R^d)$ is dense both in $(\sF, \sE_1)$ and in
$(C^\infty_\infty, \| \cdot \|_\infty)$,
 $(\sE, \sF)$ is a regular Dirichlet form on $\R^d$.
Our result now follows from \cite[Chapter 7]{FOT94}.
\qed

It is well known that Nash's inequality implies the operator
norm estimate for the transition semigroup $P_t$ from $L^1 (\R^d)$
to $L^\infty (\R^d)$. However this only implies for every $t>0$ the
existence of $p(t, x, y)$ almost everywhere on $\R^d\times \R^d$
such that for every $f\geq 0$ on $\R^d$,
$$
P_t f(x)=\int_{\R^d} p(t, x, y) f(y) dy \qquad \hbox{for a.e. }
x\in \R^d.
$$
We need something stronger.

\medskip

Since the following result has independent interest, we state
and prove it in a more  general context. For the next theorem only,
 let $E$ be a locally compact separable metric space
and $m$ a Radon measure on $E$ whose support is all of $E$. A symmetric
Dirichlet form $(\sE, \sF)$ in $L^2(E, m)$ is said to be regular if
$C_c(E)\cap \sF$ is dense both in $(\sF, \sE_1)$ and in $(C_c (E),
\| \cdot \|_\infty )$. It is well known (cf. \cite{FOT94}) that
a regular symmetric Dirichlet form $(\sE, \sF)$ has associated with it
a symmetric Hunt process $X$ that can start from every
point outside a properly exceptional set $\sN$ (cf. Theorem
\ref{thmzero}). For $x\in E\setminus \sN$, we use $\{P(t, x, dy),
t\geq 0\}$ to denote the transition probability of $X$. The
transition semigroup $\{P_t, t\geq 0\}$  of $X$ is  defined for
 $x\in E\setminus \sN $ by
$$ P_t f (x):=\E^x \left[f (X_t)\right] \qquad \hbox{for } f\geq 0
\hbox{ on $E$ and } t>0.
$$

\medskip

\begin{theorem}\label{T:qc}
Let $E$, $m$, and $P_t$ be as above.
 Assume that there is a positive left continuous
function $M(t)$ on $(0, \infty)$ such that
\begin{equation}\label{eqn:ub}
 \| P_t f\|_\infty \leq M(t) \| f\|_1  \qquad \hbox{for every } f\in L^1(E,
m) \hbox{ and } t>0.
\end{equation}
Then there is a properly exceptional set $\sN \subset E$ of $X$ and
a positive symmetric kernel $p(t,x,y)$ defined on $(0, \infty)\times
(E \setminus \sN) \times (E \setminus \sN)$
 such that $P(t,x,dy) = p(t,x,y)\,m( dy)$,
$$ p(t+s,x,y) = \int p(t,x,z) p(s,z,y) dz
\qquad \hbox{for every } x, y \in E\setminus \sN \hbox{ and } t, s
>0,
$$
  and
\begin{equation}\label{eqn:heat5}
p(t, x, y) \leq M(t) \qquad \hbox{for every } t>0 \hbox{ and } x, y
\in E\setminus \sN.
\end{equation}
Moreover, for every $t>0$ and $y\in E\setminus \sN$, $x\mapsto p(t,
x, y)$ is quasi-continuous on $E$.
\end{theorem}

\proof
 Let $\sN $ be a properly exceptional set of $X$. Recall that
 the transition semigroup $\{P_t, t\geq 0\}$ of $X$ is defined for
 $x\in E\setminus \sN $ by
$$ P_t f (x):=\E^x \left[f (X_t)\right] \qquad \hbox{for } f\geq 0
\hbox{ on $E$ and } t>0.
$$
Let $\{f_k, k\geq 1\}\subset C_c(E)\cap \sF$ be dense in both
$L^2(E, m)$ and $L^1(E, m)$. For each fixed $t>0$ and $k\geq 0$,
$P_t f_k$ is quasi-continuous on $E$. Thus for each $t>0$, there is
a ${\cal E}$-nest $\{F_n^{(t)}, n \geq 1\}$ consisting of an
increasing sequence of compact sets such that $P_tf_k$ is continuous
on each $F_n^{(t)}$ for every $k\geq 1$ (cf. \cite[Theorem
2.1.2]{FOT94}). Let ${\cal N}_t:=E \setminus \cup_{n=1}^\infty
F^{(t)}_n$, which is ${\cal E}$-polar and in particular has zero
$m$-measure.

Inequality (\ref{eqn:ub}) yields that for every $n\geq 1$,
$$ \sup_{x\in F^{(t)}_n} |P_t f_j (x) -P_t f_k (x)| \leq M(t)
 \| f_j -f_k\|_1.
 $$
Since $\{f_k, k\geq 1\}\subset C_c(E)\cap \sF$ is dense in $L^1(E,
m)$, it follows that $P_t f$ is continuous on each $F_n^{(t)}$ and
 $$\sup_{x\in E \setminus {\cal N}_t} |P_tf(x)| \leq
M(t) \| f\|_1
$$
for every $ f\in L^1(E, m)$. Therefore for every $t>0$ and $x\in E
\setminus {\cal N}_t$, there is an integrable kernel $y\mapsto
p_0(t, x, y)$ defined $m$-a.e. on $E$ such that
\begin{equation}\label{eqn:heat2}
 \E^x \left[ f(X_t)\right]=P_t f(x) = \int_{E} p_0(t, x, y) f(y) dy
\qquad \hbox{for every } f\in L^1(E, m)
\end{equation}
and
\begin{equation}\label{eqn:heat}
 p_0(t, x, y) \leq M(t) \qquad
 \hbox{for $m$-a.e. } y\in E
\end{equation}
From the semigroup property $P_{t+s}=P_tP_s$, we have for every $t,
s>0$ and $x\in E\setminus (\sN_{t+s}\cup \sN_t)$,
\begin{equation}\label{eqn:heat3}
p_0(t+s, x, y)=\int_E p_0(t, x, z) p_0(s, z, y) m(dy) \qquad
\hbox{for } m \hbox{-a.e. } y\in E.
\end{equation}
 Note that since $P_t$ is symmetric, we have for each fixed
$t>0$,
$$p_0(t, x, y)=p_0(t, y, x) \qquad \hbox{for $m$-a.e. }
 (x, y)\in E\times E .
 $$

By enlarging the properly exceptional set $\sN$ if necessary, we may
and do assume that $\sN \supset \cup_{t\in \Q_+}\sN_t$. For every
$t>0$ and $x, y \in E \setminus \sN$, let $s\in \Q_+$ be less than
$t/3$ and define
 \begin{equation}\label{eqn:newp}
  p(t, x, y):= \int_{E} p_0(s, x, w)
\left( \int_{E} p_0(t-2s, w, z)p_0(s, y, z) m(dz) \right) m(dw).
\end{equation}
By (\ref{eqn:heat3}) the above definition is independent
of the choice of $s\in \Q_+\cap (0, t/3)$.
 Clearly, $p(t, x, y)=p(t,
y, x)$ for every $x, y\in E \setminus \sN$.
 By the semigroup property and (\ref{eqn:heat2}), we have for $\phi \geq
0$ on $E$ and $x\in E\setminus \sN$,
\begin{eqnarray*}
&& \E^x \left[ \phi (X_t)\right] \\
&=& \int_{E} \left(
 \int_{E} p_0(s, x, w) \left( \int_{E} p_0(t-2s, w, z) p_0(s, z, y) m( dz )\right)
 m(dw) \right) \phi (y) m(dy) \\
 &=& \int_{E} \left(
 \int_{E} p_0(s, x, w) \left( \int_{E} p_0(t-2s, w, z) p_0(s, y, z) m(dz) \right)
 m(dw) \right) \phi(y) m(dy ) \\
&=& \int_{E} p(t, x, y) \phi(y) m(dy).
\end{eqnarray*}
Thus $p(t, x, y)$ coincides with  $p_0(t, x, y)$ $m$-a.e. on
$E\times E$. Note that it follows from (\ref{eqn:heat}) and
(\ref{eqn:newp}) that for every $t>0$ and $x, y \in E \setminus
\sN$,
$$ p(t, x, y) \leq M(t-2s)
\qquad \hbox{for every } s\in \Q_+ \hbox{ and } s<t/3.
$$
Taking $s\downarrow 0$ yields
$$ p(t, x, y) \leq M(t)
\qquad \hbox{for every } t>0 \hbox{ and } x, y \in E \setminus \sN.
$$
For $t, s>0$ and $x, y\in E\setminus \sN$, take $s_0\in \Q_+\cap (0,
\,  (t\wedge s)/3)$, and we have by (\ref{eqn:heat3})-(\ref{eqn:newp})
\begin{eqnarray*}
&& p(t+s, x, y)\\
 &=& \int_{E} p_0(s_0, x, w) \left( \int_{E}
p_0(t+s-2s_0, w, z) p_0(s_0, y, z) m(dz) \right) m(dw) \\
&=&\int_{E^5} p_0(s_0, x, w)  p_0(t-2s_0, w, u_1) p_0(s_0, u_1, u_2)
p_0(s_0, u_2, v) p_0(s-2s_0, v, z) \\
&& \ \ p_0(s_0, y, z) m(dw) m(du_1) m(du_2) m(dz) m(dv) \\
  &=& \int_E p(t, x, v) p(s, v, y) m(dv) .
\end{eqnarray*}

We may assume that there is an $\sE$-nest $\{F_n, n\geq 1\}$ such
that $\sN=E \setminus \left(\cup_{n=1}^\infty F_n \right)$ and
$P_tf_k$ is continuous on $F_n$ for each $k\geq 1$, each $t$
rational, and each $n$. It follows from inequality (\ref{eqn:ub}),
$$ \sup_{x\in F_n} |P_t f_j (x) -P_t f_k (x)| \leq M(t)
 \| f_j -f_k\|_1
 $$
for every $t\in \Q_+$ and  $n, k\geq 1$. Since $\{f_j, j\geq 1\}$ is
dense in $L^1(E, m)$, we conclude that $P_tf$ is continuous on each
$F_n$ whenever $f\in L^1(E, m)$. By (\ref{eqn:heat}), the function
$w \mapsto \int_{E} p_0(t-2s, w, z)p_0(s, y, z) dz $ is
$L^1$-integrable on $E$, and so as a function of $x$,  $p(t, x, y)$
is continuous
  on each $F_n$ for every real $t>0$ and $y\in E \setminus \sN$.
 This proves the theorem. \qed

\bigskip

In order to get off-diagonal estimates for $p(t, x, y)$ from the
on-diagonal estimate (\ref{eqn:heat5}), we need the following.

\begin{theorem}\label{T:qc2}
Let the heat kernel $p(t, x, y)$ and the properly exceptional set
$\sN$ be as in Theorem \ref{T:qc}. Suppose that $\psi \in C_c(E)$
and that there is a positive left continuous function $M_\psi(t)$ on $(0,
\infty)$ such that
\begin{equation}\label{eqn:ub2}
 \| P^\psi_t f\|_\infty \leq M_\psi(t) \| f\|_1  \qquad \hbox{for every }
f\in L^1(E, m) \hbox{ and } t>0,
\end{equation}
where $\{P^\psi_t, t\geq 0\}$ is the semigroup defined by $P^\psi
f(x): = e^{\psi(x)} P_t ( e^{-\psi} f)(x)$.
 Then
$$ p(t, x, y) \leq e^{-\psi(x) + \psi(y)} M_\psi(t) \qquad \hbox{for every
} t>0 \hbox{ and } x, y \in E\setminus \sN.
$$
\end{theorem}

\pf Clearly by Theorem \ref{T:qc}, $\{P^\psi_t, t\geq 0\}$ admits
a heat kernel
$$ p^\psi (t, x, y):= e^{\psi (x)} p(t, x, y) e^{-\psi (y)} \qquad
\hbox{for } t>0 \hbox{ and } x, y \in E\setminus \sN.
$$
Since $\psi \in C_c(E)$, for every $x\in E\setminus E$ and $s>0$,
$y\mapsto p^\psi (t, x, y)$ is $L^1$-integrable. Thus by
(\ref{eqn:ub2}), for every $s\in (0, t)$ and $x, y\in E\setminus
\sN$,
\begin{eqnarray*}
p^\psi (t, x, y)&=& \int_E p^\psi (t-s, x, z) p^\psi (s, z, y)
m(dz)\\
& \leq & M_\psi (t-s) \int_E p^\psi (s, z, y) dz \\
&=& M_\psi (t-s) e^\psi (x) \E_x \left[ e^{-\psi (X_s)} \right].
\end{eqnarray*}
Since $M(t)$ is left continuous and $\psi \in C_c(E)$, letting
$s\downarrow 0$, we have by the bounded convergence theorem that
$$ p^\psi (t, x) \leq M_\psi (t) e^{\psi (x)} e^{-\psi (x)} =M_\psi
(t),
$$
and the conclusion of the theorem follows. \qed

\bigskip

 We are now ready to prove Theorem \ref{rough-above-lem}.
 For $0<s<1$, we use $\H^s(\R^d)$ to denote
 the usual Sobolev space of fractional order:
\begin{equation}\label{fractord}
 \H^{s} (\R^d):= \left\{ v\in L^2(\R^d, dx): \,
 \il_{\R^d\times \R^d} \frac{|v(x) - v(y)|^2}{|x-y|^{d+2s}} \; dy \; dx < \infty
\right\}.
\end{equation}

\bigskip

\noindent {\bf Proof of Theorem \ref{rough-above-lem}:}
We begin with the following inequality of Nash form:
for all functions $u \in \H^\frac{\alpha}{2}(\R^d) \cap L^{1}(\R^d)$
\begin{equation} \label{sobolev}
\left(\il_{\R^d} u^2 \; dx \right)^{1+\frac{\alpha}{d}} \leq c_1
  \left(\il_{\R^d} \il_{\R^d} \frac{|u(x) - u(y)|^2}{|x-y|^{d+\alpha}}
          \; dy \; dx \right) \;
  \left(\il_{\R^d} |u (x)| \; dx \right)^{\frac{2\alpha}{d}} \;,
\end{equation}
where $c_1$ is a positive constant depending only on the space dimension $d$.
This may be proved using the continuous Sobolev embedding
$\H^\frac{\alpha}{2}(\R^d) \hookrightarrow L^{2d/(d-\frac{\alpha}{2})}(\R^d)$
and interpolation in $L^p(\R^d)$ spaces.
An alternative way of proving this is to recall that the transition
densities for a symmetric stable process of order $\alpha$ are bounded
by $ct^{-d/\alpha}$ and then to apply Theorem 3.25 of \cite{CKS87}.

We then deduce from (A3)
\begin{eqnarray}
\| u\|_2^{2+(2\alpha/d)}
&\leq & c_1 \left( \il_{|x-y|<1} \frac{|u(x) - u(y)|^2}{|x-y|^{d+\alpha} }
          \; dy \; dx + c_2 \int_{\R^d} u(x)^2 dx \right) \| u\|_1^{2\alpha/d}
          \nonumber \\
&\leq & c_1\big(\kappa_1^{-1} \sE(u,u)+c_2 \|u\|_2^2\big) \|
u\|_1^{2\alpha/d}. \label{eqn:2.2a}
\end{eqnarray}
Let $\{P_t, t\geq 0\}$ denote the transition semigroup of $X$; that
is,
$$ P_t f(x):=\E^x \left[f (X_t)\right] \qquad \hbox{for } f\geq 0
\hbox{ on $\R^d$ and } t>0.
$$
It follows from  Theorem 2.1 of \cite{CKS87} that
 \begin{equation}\label{eqn:2.2}
  e^{-c_2 \kappa_1
t} \| P_t f\|_\infty \leq c_3 \, t^{-d/\alpha} \| f\|_1.
\end{equation}
Noting that $(\sE, \sF)$ is a regular Dirichlet form on $\R^d$, the
conclusion of this theorem now follows immediately from Theorem
\ref{T:qc}. \qed

\bigskip

Let $B\subset \R^d$ be a
ball. Denote by $X^B$ be the subprocess of $X$ killed upon leaving
$B$. Let $\{P^B (t, x, dy), t> 0\}$ be the transition probability of
$X^B$.
We will need the existence and regularity of the
transition density of $X^B$.

\begin{theorem}\label{T:2.1} Assume conditions (A1)--(A3) hold.
Let $\sN$ be the properly exceptional set of $X$ in Theorem
\ref{rough-above-lem}. There exist a positive
    symmetric kernel
$p^B(t,x,y)$ defined on $(0, \infty)\times (B \setminus \sN) \times
(B \setminus \sN)$
 such that $P^B(t,x,dy) = p^B(t,x,y) dy$, and
\begin{equation}
\label{pbub}
p^B(t, x, y) \leq C_1 t^{-d/\alpha} \qquad \hbox{for every } t>0
\hbox{ and } x, y \in B\setminus \sN \,,
\end{equation}
where the constant $C_1$
depends on $\alpha, \beta, \kappa_1, d$. Moreover, for every $t>0$
and $y\in B\setminus \sN$, $x\mapsto p^B(t, x, y)$ is
quasi-continuous on $B$.
\end{theorem}

\pf Let $\sN$  and $p(t, x, y)$  be the properly exceptional set and
the transition density function, resp., in Theorem \ref{rough-above-lem}.
Define $\tau_B=\inf\{t>0: X_t\notin B\}$. Then
$$ p^B(t, x, y):= p(t, x, y) -\E^x \left[ p(t-\tau_B, X_{\tau_B}, y);
\tau_B<t \right], \qquad x, y \in B\setminus \sN,
$$
is the transition density function for $X^B$. It is easy to see that
$p^B (t, x, y)$ is symmetric and $y\mapsto p(t, x, y)$ is
quasi-continuous.

Note that the Dirichlet form for $X^B$ is $(\sE, \sF^B)$, where
\begin{equation}\label{eqn:c3}
 \sF^B =\{u\in \sF: u=0 \hbox{ $\sE$-q.e. on } B^c \}.
\end{equation}
So for $u\in \sF^B$,
\begin{equation}\label{eqn:c4}
 \sE(u, u) = \int_{B\times B} (u(x)-u(y))^2 J(x, y) dx dy
 + \int_B u(x)^2 \kappa_B(x) dx,
 \end{equation}
 where \begin{equation}\label{eqn:c5}
\kappa_B (x) = 2 \int_{B^c} J(x, y) dy.
\end{equation}
 It follows from
 (A3) that there exists a constant $c_1>0$ such that
   $$c_1\dist(x, \partial B)^{-\al}
   \leq \kappa_B(x) \leq c_1^{-1} \, {\rm dist}(x, \partial B
)^{-\beta}
  \qquad \hbox{for } x \in B.
 $$
Thus we have from (\ref{eqn:2.2a}) that there is a constant $c_2>0$
such that for $u\in \sF^B$,
$$ \| u\|_2^{2+(2\alpha/d)} \leq c_2 \sE(u, u) \, \| u\|_1^{2\alpha /d}.
$$
It follows from Theorem 2.1 of [CKS87] that the transition semigroup
$\{P^B_t, t>0\}$ of $X^B$ satisfies
$$ \| P^B_t f \|_\infty \leq c_3 t^{-d/\alpha} \| f \|_1 \qquad
\hbox{for } f\in L^1(B, dx) \hbox{ and } t>0.
$$
This implies that
$$ p^B(t, x, y) \leq c_3 t^{-d/\alpha} \qquad \hbox{for } t>0
\hbox{ and } x, y \in B\setminus \sN.
$$
\qed

\bigskip

\begin{remark}\label{meyer}
{\rm We will use several times the following construction of Meyer
\cite{Mey75} for jump processes. Suppose we have a jump intensity
kernel $J(x,y)$  and another jump intensity kernel $J_0(x,y)\leq
J(x,y)$ such that
\begin{equation}\label{sjdef}
 \sJ(x):=\int_{\R^d} (J(x, y)-J_0(x,y))\, dy \le c_1   \q \hbox{for all $x$}.
\end{equation}
Let
\begin{equation}\label{qdef}
  q(x,y) = \frac{J(x, y)-J_0( x,y)}{\sJ(x)}.
\end{equation}
Let $Z^{(0)}=\{Z^{(0)}_t, t\geq 0\}$ be the process corresponding to
the jump kernel $J_0$. Then we can construct
a process $Z$ corresponding to the jump kernel $J$ as follows.
Let $S_1$ be an  exponential random variable of parameter 1 independent
   of $Z^0$,
let $C_t=\int _0^t \sJ(Z^{(0)}_s)\, ds$, and let $U_1$ be the first
time that $C_t$ exceeds $S_1$. We let $Z_s = Z^{(0)}_s$ for $0\le s
\le U_1$.

At time $U_1$ we introduce a jump from $Z_{U_1-}$ to $Y_1$, where
$Y_1$ is chosen at random according to the distribution $q(Z_{U_1-},\,y)\,dy$.
We set $Z_{U_1}=Y_1$, and repeat, using an independent
exponential $S_2$, etc.  Since $\sJ(x)$ is bounded, only finitely many
new jumps are introduced in any bounded time interval.  In
\cite{Mey75} it is proved that the resulting process corresponds to
the kernel $J$. See also \cite{INW66}.

Note that if $\sN_0$ is the null set corresponding to $Z^{(0)}$ then
this construction yields that $\sN \subset \sN_0$.
}\end{remark}

\medskip
\begin{remark}\label{conversemeyer}
{\rm
Conversely, we can also remove a finite number of jumps from a process
$Z$ to obtain a new process $Z^{(0)}$.  Suppose
$J(x,y)=J_0(x,y)+J_1(x,y)$, where we have $\int J_1(x,y)\, dy\leq c_1$
for all $x$ and for simplicity we also assume that
$J_0(x, y) J_1 (x,y)=0$.
One starts with the process $Z$ (associated with $J$), runs it
until the stopping time $S_1=\inf\{t: J_1(Z_{t-},Z_t)>0\}$, and at
that time restarts $Z$ at the point $Z_{S_1-}$. One then
repeats this procedure.
Meyer \cite{Mey75} proves that
the resulting process $Z^{(0)}$ will correspond to the jump kernel
$J_0$.  In this case we have $\sN_0 \subset \sN$. }
\end{remark}

We will need the following bound, which arises easily from Remark \ref{meyer}.

\begin{lemma} \label{mey-bound}
Let $Z^{(0)}$, $Z$ be as in Remark \ref{meyer}, and
$F\in \sigma( Z^{(0)}_t, 0\le t <\infty)$. Then
$$ \P^x( \{Z_s = Z^{(0)}_s \hbox{ for all } 0\le s \le t\} \cap F) \ge
e^{- t ||\sJ||_\infty}\bP^x(F). $$
In particular, if $D \subset \bR^d$ and $A \subset D$,
\begin{equation}\label{meybineq}
 \P^x( Z_t \in A, \tau^Z_D > t ) \ge
e^{- t ||\sJ||_\infty}\P^x( Z^{(0)}_t \in A, \tau^{Z^{(0)}}_D > t ).
\end{equation}
\end{lemma}

\proof We have, writing $c_1=||\sJ||_\infty$ and
$G= \{Z_s = Z^{(0)}_s, 0\le s\le t\}$,
$$ \bP^x(G \cap F)= \P^x( U_1 >t,F )= \P^x( C_t < S_1,F)
\ge \P^x(S_1 >  c_1 t,F)=\bP^x(F) e^{-c_1t}. $$
The last equality holds because $S_1$ is independent of the process $Z^{(0)}$.

For the second assertion, let $F=\{ Z^{(0)}_t \in A,  \tau^{Z^{(0)}}_D > t\}$,
and $G=\{ Z_s = Z^{(0)}_s \hbox{ for all } 0\le s \le t\}$. Then,
$$  \P^x( Z_t \in A, \tau^Z_D > t )
 \ge  \P^x( Z_t \in A, \tau^Z_D > t, G )= \bP^x(F \cap G)
\ge e^{- t ||\sJ||_\infty} \bP^x(F). $$
\qed

Let $J(x,y)=J_0(x,y)+J_1(x,y)$, and let $\sJ$, $q$ be defined
by (\ref{sjdef}), (\ref{qdef}). Write $p^{(0)}(t,x,y)$ for the
transition density of the process $Z^{(0)}$ associated with $J_0$.
Let
\begin{equation}
 r  (t,x,y) = \int q(x,z) p(t,z,y) dz.
\end{equation}

The following lemma, which follows quite easily from Meyer's
construction, is proved in \cite{BGK}.

\begin{lemma}\label{jumpf}
(a) For any Borel set $B$
\begin{equation}\label{cjmp}
  \bP^x( Z_t \in B) = \bP^x(   Z^{(0)}_t \in B, S_1>t)
 + \bE^x \int_0^t \int_B r(t-s,Z^{(0)}_s,z) \sJ(Z^{(0)}_s)\, dz \,ds.
\end{equation}
(b) If $\| J_1 \|_\infty < \infty $ then
\begin{equation}\label{ptin}
  p(t,x,y) \le p^{(0)}(t,x,y) + t \| J_1 \|_\infty.
\end{equation}
\end{lemma}

We need the following standard result.

\begin{lemma}\label{there-back}
Suppose that there exist positive constants $r$, $t$ and $p$ such
that
\begin{equation}\label{tba}
 \bP^x( |X_s -x| > r) \le p \qquad \hbox{ for all } x \in \bR^d\setminus \sN
  \ \hbox{ and } \ 0\le s\le t.
\end{equation}
Then
\begin{equation*}
   \bP^x( \sup_{0\le s \le t} |X_s -X_0| > 2r )
   \le 2p \q \hbox{ for all } x \in \bR^d\setminus \sN.
\end{equation*}
\end{lemma}

\proof Let
  $S=\inf\{t: |X_t-X_0|> 2r\}$.
Then using the strong Markov property of $X$ and (\ref{tba})
\begin{align*}
\bP^x \Big(\sup_{s\le t}|X_s-x| > 2r \Big) &= \bP^x(S \le t)\\
&\le \P^x(|X_{t}-x|>r) + \P^x(S\leq t \hbox{ and }  X_{t}\in B(x,r)) \\
&\leq p + \bP^x(S\le t \hbox{ and }  |X_{t}-X_S| > r) \\
&= p + \bE^x \left( \mathbbm{1}_{\{ S \le t\}}
\bP^{X_S}(|X_{t-S}-X_0| > r) \right) \\
& \le 2p.
\end{align*}
\qed

We now use Lemma \ref{jumpf} to obtain off-diagonal upper bounds
on $p(t,x,y)$; this technique was introduced in \cite{BGK}.

\begin{proposition} \label{exittime}
Let $\sN$ be the properly exceptional set of $X$ in Theorem
\ref{rough-above-lem}. There exist
   constants  $t_0>0$, $c_1, c_2,
c_3>0$ such that for every $x,y \in \R^d \setminus \sN$,
\begin{description}
\item{\rm (a)} \  $p(t, x, y)\leq c_1$ if
 $|x-y|\geq 1$ and $t\le 2$.
\item{\rm (b)}\  \
$\displaystyle \P^x \left(\sup_{s\leq t} |X_s-x|> r \right)< c_2 t
e^{-c_3r}$
 for every $r \ge 1/8$, $t\le t_0$.
\item{\rm (c)} There exists $t_1>0$ such that
\begin{equation}\label{neededlb}
 \bP^x \left(\sup_{s\leq t_1} |X_s-x|> 1/4 \right)<1/4.
\end{equation}
 \end{description}
\end{proposition}

\proof (a) For $\delta>0$ let
$J^{(\delta)}(x,y)=J(x,y) \mathbbm{1}_{(|x-y|<\delta)}$, and let $\sE^{(\delta)}$ be
defined by (\ref{Dirform}) with $J^{(\delta)}$ in place of $J$.
(Of course, if $\delta>1$ then  $\sE^{(\delta)}$ is just $\sE$.)
Clearly $(\sE^{(\delta)}, \sF)$ is a regular Dirichlet form on
$\R^d$. Thus there is a Hunt process $X^{(\delta)}$ associated with
$(\sE^{(\delta)},\sF)$ and a properly exceptional set $\sN(\delta)$
so that $X$ starts from every point in $\R^d
    \setminus \sN(\delta)$.
Since using Meyer's procedure we can construct  $X^{(\delta)}$ from
$X$ and vice versa, we can take $\sN(\delta)=\sN$. Let $p^{(\delta)}
(t, x, y)$ be the  transition density of $X^{(\delta)}$.

By (A1)--(A3) we have
\begin{eqnarray*}\label{jdbounds}
  \int_{B(x,\delta)^c} J(x,y)\, dy
&=& \int (J(x,y)- J^{(\delta)}(x,y))dy  \le c_4 \delta^{-\beta}, \\
  \int_{B(x,\delta)} |x-y|^2 J(x,y) dy  &=& \int  |x-y|^2
   J^{(\delta)}(x, y) \, dy
\le c_5 \delta^{2-\beta}. \\
 || J -  J^{(\delta)}||_\infty
   &=& \sup_{ |x-y|\ge \delta } J(x,y) \le c_6 \delta^{-d-\beta}. \\
  \end{eqnarray*}

Starting from (\ref{sobolev}) we have
\begin{align*}
\| u\|_2^{2+(2\alpha/d)}
&\leq c_7 \Big( \int \int_{|x-y|<\delta} \frac{|u(x)-u(y)|^2}{|x-y|^{d+\alpha}}
\, dy\, dx+\delta^{-\alpha} \| u\|_2^2\Big)\| u\|_1^{2\alpha/d}\\
&\leq c_8(\kappa_1^{-1} \sE^{(\delta)}(u,u)
   +\delta^{-\al}\| u\|_2^2) \|u \|_1^{2\alpha/d}.
\end{align*}
Theorem \ref{T:qc2} in Section 2 and Theorem 3.25 of \cite{CKS87}
now give
\begin{equation}\label{ex1}
p^{(\delta)}(t,x,y)\leq c_9 t^{-d/\alpha} e^{c_{10} t\delta^{-\alpha}}e^{-E_\delta(2t,x,y)}
\qquad \hbox{for } x, y \in \R^d \setminus \sN(\delta),
\end{equation}
where
\begin{align*}
\Gamma_\delta(f,f)(x)&= \int (f(x)-f(y))^2 J^{(\delta)}(x,y)\, dy,\\
\Lambda_\delta(\psi)^2&=\| e^{-2\psi} \Gamma_\delta(e^\psi, e^\psi)\|_\infty
\lor \| e^{2\psi} \Gamma_\delta(e^{-\psi},e^{-\psi})\|_\infty,\\
E_\delta(t,x,y)&=\sup\left\{|\psi(x)-\psi(y)|-t\Lambda_\delta(\psi)^2:
\ \psi \in \sF \cap C_b(\R^d) \hbox{ with } \Lambda_\delta(\psi)<\infty \right\}.
\end{align*}

Let $t>0$ and $x_0, y_0 \in \R^d \setminus \sN$; write $R=|x_0-y_0|$.
Note that the set $\sC$ of compactly supported Lipschitz continuous functions
is a core for $\sE^{(\delta)}$.
Let $\lam>0$ and
$$ \psi(x) = \lam(R-|x_0-x|)^+. $$
So $|\psi(x)-\psi(y)| \le \lam |x-y|$.  Noting that
$|e^t-1|^2\le t^2 e^{2|t|}$, we obtain
\begin{eqnarray*}
\Lambda_\delta(e^\psi)^2(x) &=& \int(e^{\psi(x)-\psi(y)}-1)^2 J^{(\delta)}(x,y)dy \\
 &\le & e^{2\lam \delta} \lam^2 \int |x-y|^2  J^{(\delta)}(x,y)dy  \\
 &\le & c_{11} (\lam \delta)^2  e^{2\lam \delta} \delta^{-\beta}
  \le  c_{12} e^{3\lam \delta} \delta^{-\beta}.
  \end{eqnarray*}
Since the same bound holds for $\Lambda_\delta(e^{-\psi})^2(x)$ we have
\begin{equation}
- E_\delta(2t,x_0,y_0) \le -\lam R + c_{12} t  \delta^{-\beta}
e^{3\lam \delta}.
 \end{equation}

In what follows we will always choose $\delta$ so that $\delta\le e$, and $t$
such that $t<\delta^\beta$. Hence $t \le \delta^\beta\le c_{13} \delta^\al$, so that
the term $e^{c_{10} t\delta^{-\alpha}}$ in (\ref{ex1}) is bounded by $c_{14}$.
We take
\begin{equation}
 \lam = \frac{1}{3\delta} \log \left(\frac{\delta^\beta}{t}\right)>0.
\end{equation}
So
\begin{equation}\label{ptubbase}
  - E_\delta(2t,x_0,y_0) \le  - \frac{R}{3\delta} \log \Big(\frac{\delta^\beta}{t} \Big)
 + c_{15} t \delta^{-\beta} \Big( \frac{\delta^\beta}{t}\Big)
 =  \frac{R}{3\delta} \log \Big(\frac{t}{\delta^\beta}\Big) + c_{16} .
\end{equation}
Hence by (\ref{ex1}),
\begin{equation} \label{ptuba}
   p^{(\delta)}(t,x_0,y_0)
 \le  c_{17} t^{-d/\al}  \Big(\frac{t}{\delta^\beta}\Big)^{R/(3\delta)}.
 \end{equation}

We need to consider two cases. Set $R_0=3e (1+d/\al)$.
First, suppose $R\ge R_0$ and $t\le 1$. Set $\delta=e$. Then
since there are no jumps larger than 1, we have
$p(t,x_0,y_0)=p^{(\delta)}(t,x_0,y_0)$, and by (\ref{ptuba})
\begin{equation} \label{ex4}
 p(t,x_0,y_0) \le c_{17} t^{-d/\al}  \Big(\frac{t}{\delta^\beta} \Big)^{R/(3\delta)}
    \le c_{18} t e^{- \beta R/3e}.
\end{equation}
Next, suppose  $0<R\le R_0$, and set
$\delta= R\al/(3(\al+d))$. We assume that
$ t < \delta^\beta = c_{19}R^\beta$.
Then since $R/(3\delta)= 1+d/\al$, (\ref{ptin}) and \ref{ptuba})  give
\begin{equation}\label{ptubc}
 p(t,x_0,y_0) \le c_{20} \frac{t}{\delta^{\beta+ \beta d/\al}}
 +  c_{21} \frac{ t}{  \delta^{\beta+d}}
 \le c_{22} \frac{t}{R^{\beta+ \beta d/\al}}
 +  c_{23} \frac{ t}{  R^{\beta+d}}
\end{equation}
Combining (\ref{ex4}) and (\ref{ptubc}), we deduce there exists
$t_1>0$ such that if $|x_0-y_0|=R$ and
 $R\ge 1/16$, and $t\le t_1$,
   then
\begin{equation} \label{ptubb}
 p(t,x_0,y_0) \le c_{24} t e^{-c_{25} R}  \qquad
 \hbox{for }  t\in (0,  t_1 ].
\end{equation}
(a) now follows on combining (\ref{ptubb}) and (\ref{rough-above}).

\medskip Integrating (\ref{ptubb}) and (\ref{ex4}) over
$B(x, \, R)^c$, we obtain, for $x\not\in \sN$,
\begin{equation*}
  \bP^x( X_t \not\in B(x, R)) \le c_{26} t e^{-c_{27}R}
\qquad  \hbox{ if } R \ge 1/16 \hbox{ and }  t \in (0,  t_1 ],
\end{equation*}
and (b) follows by Lemma \ref{there-back}.
(c) is immediate from (b). \qed

\begin{theorem}\label{reg_conserv}  The process $X$ is conservative.
\end{theorem}

\noindent{\bf Proof:} From the probabilistic interpretation of the
Beurling-Deny decomposition for $(\sE, \sF)$ (cf. \cite{FOT94}),
$X$ admits no killing inside $\R^d$.
On the other hand, by Proposition \ref{exittime} the process $X$ cannot
escape to infinity in the time interval $[0,1]$. Since the lifetime $\zeta$
of $X$ satisfies $\bP^x(\zeta\ge 1)=1$ for all $x\in \bR^d\setminus \sN$,
we have  $\bP^x(\zeta =\infty)=1$, hence $X$ is conservative.
\qed

\section{Lower bounds on the heat kernel}

Throughout this section we will assume that the jump kernel $J$
satisfies (A1)--(A3) and also \a4xi for some $\xi\in (0,1)$, unless
otherwise specified.

\begin{lemma}\label{init-m-bounds} Let $J(x,y)$ satisfy
(A1)--(A3), and \a4xi for some $\xi\in (0,1)$.
Then
$$ \sF=\{ f\in L^2(\R^d, dx): \, \sE (f, f)<\infty\}
= \H^{\beta/2}(\R^d).  $$
The corresponding Hunt process $X$ is a
conservative strong Markov process which can start from any point in
$\R^d$.
\end{lemma}

\proof For two bilinear forms $\sE$
and $\sC$ having a common core $C^1_c(\R^d)$, let us write
$\sE (f, f) \approx \sC (f, f)$  if
there is a finite constant $c_1>0$ such that
$c_1^{-1} \sE(f, f) \leq \sC (f, f) \leq c_1 \sE (f, f)$ for every $f\in C^1_c(\R^d)$.
As $(\sE, \sF)$ satisfies (A1)--(A3) and \a4xi, we
have for $f\in C^1_c(\R^d)$,
\begin{eqnarray*}
\sE_1 (f, f) &\approx& \int_{|x-y|\leq \xi} \frac{(f(x)-f(y))^2}{|x-y|^{d+\beta}} dx\,dy
   + \| f\|_2^2 \\
&\approx& \int_{\R^d\times \R^d} \frac{(f(x)-f(y))^2}{|x-y|^{d+\beta}} dx\,dy
   + \| f\|_2^2 \\
&:=& \sC_1 (f, f).
\end{eqnarray*}
So $\sF=\ol{C^1_c(\R^d)}^{\sE_1}=\ol{C^1_c(\R^d)}^{\sC_1}=\H^{\beta/2}(\R^d)
=\{f: \sC_1(f, f)<\infty \}$,
which is the same as
$ \{f: \sE_1(f, f)<\infty \}$.

Note that the process $X$ can be constructed from
 the L\'evy process $Z$ on $\R^d$ whose L\'evy measure is
$$ J_0(h) \, dh=|h|^{-d-\beta}\mathbbm{1}_{(|h|\leq \xi)}\,dh   $$
using Remark \ref{meyer}.
As $Z$ is conservative  and can start from every point in $\R^d$,
the same is true of $X$. \qed

\medskip
In view of the above we can take the properly exceptional
set $\sN$ to be the empty set throughout this section.

\smallskip
Let $B$ be a ball of radius $R$ centered at 0 for $R\in [1,4]$.
Let $p^B(t, x, y)$ be the transition density function of the subprocess
$X^B$ of $X$ killed upon exiting the ball $B$.

\begin{lemma}\label{delut}
There exists a positive constant $c_1$ depending on $\xi$ such that
$$ p^B(t, x, y)\leq c_1 t^{-d/\alpha}
\quad \hbox{ and } \quad
\left|\frac{\partial p^B(t, x, y)}{\partial t}
\right|\leq c_1 t^{-1-\frac{d}{\alpha}}  $$
for every $x, y\in B$ and $t>0$.
\end{lemma}

\proof
The assertion for $p^B$ follows immediately  from (\ref{rough-above}). As
$$\int_{B\times B} p^B(t, x, y)^2 dx\,
dy = \int_B p^B(2t, x, x) dx<\infty, $$
the symmetric semigroup
$P^B_t$ of $X^B$ is a Hilbert-Schmidt operator on $L^2(B, dx)$ and
so it is compact and has discrete spectrum $\{e^{-\lambda_i t},
i\geq 1\}$, with repetitions  according to multiplicity. Let
$\{\psi_i, i\geq 1\}$ be the corresponding eigenfunctions normalized
to have unit $L^2$-norm on $B$ and to be orthogonal to each other.
Then it is well known (see, e.g., \cite{Bas97}, Section VII.6) that
$$ p^B(t,x, y)=  \sum_{i=1}^\infty e^{-\lam_i t} \psi_i(x)\psi_i(y). $$
Hence
\begin{equation}\label{eqn:c2}
\frac{\partial p^B(t, x, y)}{\partial t}(t,x)
=-\sum \lam_i e^{-\lam_i t} \psi_i(x)\psi_i(y),
\end{equation}
where the convergence is in the $L^2$--sense.
Note that for any given $\delta >0$,
by Cauchy-Schwarz and the Chapman-Kolmogorov equations,
\begin{align*}
 \psi_i (x) &= e^{\lambda_i \delta} \int_B p^B(\delta, x, y) \psi_i (y) dy\\
&\leq e^{\lam_i \delta} \Big(\int_B p^B(\delta,x,y)^2\, dy\Big)^{1/2}
\Big(\int_B (\psi_i(y))^2\, dy\Big)^{1/2}\\
  &\leq e^{\lambda_i \delta} p(2\delta, x, x)^{1/2}
\leq c_2 \delta^{-d/\alpha} e^{(c_3+ \lambda_i) \delta}
\end{align*}
for every $x\in B$ by (\ref{rough-above}).
It follows that the series on the right hand side of
(\ref{eqn:c2}) converges uniformly on
$(2\delta, \delta^{-1})\times B \times B$.
The function $x\to xe^{-xt/2}$ is bounded over nonnegative $x$ by $c_4/t$, so
for $x, y \in B$, the
right hand side of (\ref{eqn:c2}) is bounded in absolute value by
\begin{align*}
\frac{c_4}{t} \sum_i e^{-\lam_i t/2}|\psi_i(x)|\, |\psi_i(y)|
&\leq \frac{c_4}{t}\Big( \sum_i e^{-\lam_i t/2} \psi_i(x)^2\Big)^{1/2}
\Big(\sum_i e^{-\lam_i t/2} \psi_i( y)^2\Big)^{1/2}\\
&= \frac{c_3}{t} p^B(t/2,x,x)^{1/2} p^B(t/2, y, y)^{1/2}.
\end{align*}
Our upper bound (\ref{pbub}) on $p^B(t, \cdot, \cdot)$ yields the
desired bound on  $|\partial p^B(t, x, y)/\partial t|.$
\qed

\begin{lemma} \label{stabproc}
Let $Y$ be a pure jump symmetric process on $\R^d$ with jump kernel
 $J_Y(x,y)$ satisfying
\begin{equation} \label{ybeta}
c_1 |x-y|^{-d-\beta} \le J_Y(x,y) \le c_2  |x-y|^{-d-\beta},
\qq x,y\in \R^d.
\end{equation}
Let $D=B(0,\half)$ and $p^D_Y(t,x,y)$ be the transition density of
  the subprocess of $Y$ killed upon exiting $D$.
 Then there exist constants
$c_3$, $c_4$ (depending on $c_1$, $c_2$) such that
\[  p^D_Y(t,x,y) \ge c_3 e^{-c_4 t}  \qq \text{ for }
t \ge \tfrac34 \hbox{ and } x, y \in B(0, \tfrac14). \]
\end{lemma}

\proof
It is known from \cite{CK03} that such a process $Y$ is a Feller process that
can start from every point in $\R^d$ and has a H\"older continuous
transition density function $p(t, x, y)$. This implies in particular
that the transition density function $p^D_Y(t, x, y)$ for the
subprocess of $Y$ killed upon leaving $D$ exists and is continuous
on $(0, \infty)\times D\times D$. For each fixed $x\in D$, $(t,
y)\mapsto p^D_Y(t, x, y)$ is a caloric function of $Y$ in $(0,
\infty) \times D$.

Let $B'=B(0,\tfrac14)$. Then
 $\P^0(\tau^Y_{B'} \ge \tfrac12) \ge c_5>0$. By \cite{CK03} a
parabolic Harnack inequality holds for $Y$. Therefore
  for every $x\in D$, $\tfrac14 \le t_1
   \le \tfrac12$ and $\tfrac34 \le t_2 \le 1$,
$$ p^D_Y( t_1,x,y) \le
   c_6 p^D_Y(t_2,x,y') \qq  \hbox{for } y, y' \in B'. $$
So if $x_1, y_1 \in B'$,
\begin{align*}
 c_5 &\le \int_{B'}  p^{B'}_Y(t_1,x_1,y) dy \\
     &\le \int_{B'}  p^D_Y(t_1,x_1,y) dy \\
     &\le c_6 \int_{B'}  p^D_Y(t_2,x_1,y_1) dy
     = c_6 |B'|   p^D_Y(t_2,x_1,y_1).
\end{align*}
This proves the result for $t \in [\tfrac34,1]$. An easy
iteration argument now gives the result for $t>1$. \qed

\begin{proposition}\label{utdecay}
Let $B=B(0,R)$ with $R\in [1, 4]$
 and $0<t_0<t_1 < \infty$. There
is a constant $c_1=c_1(\xi, t_0,t_1)>0$ such that
$$ p^B(t, x, y)\geq c_1(R-|x| )^{\beta}(R-|y|)^{\beta} \q
\text{ for every $ t\in [t_0,t_1]$ and  $x, y \in B$.} $$
\end{proposition}

\proof
Recall we are working under \a4xi; the constants $c_i$ in this proof will
depend on $\xi$. We can assume without loss of generality  that $\xi \le \tfrac18$.
By considering the process $X'_t = X_{2 t/t_0}$, which satisfies (A1)--(A3)
and \a4xi (but with different constants $\kappa_i)$, we can assume
that $t_0 \ge 2$.

Let $V$ be a L\'evy process with jump kernel
\begin{equation*}
  J_0(x,y) =\
\begin{cases}
 \kappa_2 |y-x|^{-d-\beta} &\hbox{ if } |y-x|\leq \xi, \\
 \kappa_1 |y-x|^{-d-\alpha} &\hbox{ if } \xi< |y-x| \le 1, \\
 0 &\hbox{ if } |y-x|>1.
\end{cases}
\end{equation*}
We can construct $X$ from $V$ by using the construction of  Remark \ref{meyer}.
Since $\sJ(x):=\int_{\R^n} (J(x, y)-J_0(x, y)) dy $ is bounded, by
Lemma \ref{mey-bound} there is a positive constant $c_2=c_2(t_1)>0$ such that
$$ \P^x \left(X_t\in A \hbox{ and } \tau_B^X>t \right) \geq c_2 \,
\P^x \left(V_t\in A \hbox{ and }  \tau_B^{V}>t \right) $$
for every  $t\in (0, t_1] $ and $A\subset B$.
It thus suffices to get a lower bound on $v^B(t, x, y)$, the transition
density of the subprocess of $V$ killed upon exiting $B$.

 By the Chapman-Kolmogorov equations,
$$v^B(t,x,y)\geq \int_{B(0,1/4)} v^B(t/2,x,z)v^B(t/2,z,y)\, dz.$$
Therefore it is enough to prove that
\begin{equation}
\label{vlowerb}
 v^B(t,y,z) \ge c_3 (R-|y|)^\beta \q
\text { for $y \in B$, $z \in B(0,1/4)$, $t\in [1,t_1]$}.
\end{equation}
  Let $y\in B$,  $\eps=R-|y|$,
   and $\eta= (\eps \wedge \xi)/2$.
Let $y_0\in B(0, R - \half)$ with $|y_0-y| \le \tfrac34$. Our first
estimate is to prove that there exist $c_4,c_5>0$ such that
\begin{equation}\label{P34.1}
\P^y \left( V_{t \eta^\beta}\in B(y_0, 1/4) \hbox{ and }  \tau^V_B
>t \eta^\beta \right) \ge c_5 t  \eta^\beta  \qq \hbox{ for } t \le c_4.
\end{equation}
To prove this, define
$$
  V_t^\eta = V_0 +\sum_{0<s\leq t} \Delta V_s \mathbbm{1}_{(|\Delta V_s|\le
  \eta)}.
$$
Note that $V^\eta$ is a L\'evy process with jump kernel
$\kappa_2 |x-w|^{-d-\beta}\mathbbm{1}_{( |x-w| \le \eta )}$
 and that
the processes $V^\eta$ and $V-V^\eta$ are independent. We write
$J_1(x,w)$ for the jump kernel of $V-V^\eta$.

 Each coordinate of
$V^{\eta}$ is a square integrable martingale. Applying Doob's
maximal inequality to each coordinate of $V^{\eta}$ separately
yields
\begin{align}
\label{vsmineq}
   \P \left(\sup_{s\leq  t \eta^\beta}
  |V^{\eta}_s-V^{\eta}_0| \ge \eta \right)
 &\le  4d \, (\eta/d)^{-2} \, \bE \left[ |V^{\eta}_{t \eta^\beta}
   -V^\eta_0|^2
  \right]  \nonumber \\
  &= 4d^3  t \eta^\beta \eta^{-2} \,
    \int_{ |h|\le \eta} \kappa_2 |h|^{2-d-\beta}
    dh \nonumber \\
    & \leq  c_6 t .
 \end{align}
Let $c_4=\min \left\{ \tfrac1{4c_6}, \, \tfrac1{24} \right\}$.

There are constants $c_8>c_7>0$ such that the total L\'evy measure
of $V-V^\eta$ is bounded by
$$  c_7 \eta^{-\beta} \leq \int_{\R^d} J_1(0,x) dx \leq c_8
\eta^{-\beta} .
$$
Since $|y_0-y| \le \tfrac34$,
\[ \int_{B(y_0-y,1/8)} J_1(0,w) dw \ge c_9 . \]
Let $t \le c_4$ and $F(t)$
be the complement of the the event on the left side of
(\ref{vsmineq}). Then $\P(F(t))\geq 3/4$.
   Let $G(t)$ be the event
that $V-V^\eta$ makes exactly one jump in the time interval $[0,t
\eta^\beta]$, and that
the size falls within the ball $B(y_0-y,1/8)$. Then
\[ \P^y(G(t))
 \ge ( c_7 \eta^{-\beta})( t\eta^\beta  ) e^{- (c_8 \eta^{-\beta})( t \eta^\beta ) }
\frac{c_9}{c_8\eta^{-\beta}}
 \ge c_{10}  t \eta^\beta . \]
 Note that $\eta \leq \xi/2\leq 1/16$.
If both $F(t)$ and $G(t)$ occur, then $V_{t \eta^\beta} \in B(y_0,
1/4)$ and $\tau^V_B > t \eta^\beta$. As $V^{\eta}$ and $V-V^{\eta}$
are independent, we obtain  (\ref{P34.1}).
In particular, we have for every $y\in B$, with $s_1:=c_4 2^{-\beta}
\left( (R-|y|)\wedge \xi \right)^\beta$,
 \begin{eqnarray*}
 \P^y \left(V_{s_1}  \in B(0, R-\tfrac14) \hbox{ and } \tau^V_B
> s_1 \right) &\geq&  \P^y \left(V_{s_1} \in
B(y_0, 1/4) \hbox{ and } \tau^V_B > s_1 \right) \\
&\geq & c_5 \, s_1.
\end{eqnarray*}
Similarly, for every $z\in B(0, R-\tfrac14 )$, let $z_0\in B(0,
R-\tfrac{3}{4})$ with $|z-z_0|<3/4$. Note that in this case
$R-|z|\geq 1/4$ and $\xi \leq 1/8$ and so $\eta:=\frac{(R-|z|)\wedge
\xi}2=\tfrac{\xi}2$.
 The same argument as above shows
that, with $s_2:=c_4 2^{-\beta}\xi^{\beta}$,
 \begin{eqnarray*}
 \P^z \left(V_{s_2}  \in B(0, R-\tfrac12) \hbox{ and } \tau^V_B
> s_2 \right) &\geq&  \P^z \left(V_{s_2} \in
B(z_0, 1/4) \hbox{ and } \tau^V_B > s_2 \right) \\
&\geq & c_5 \, s_2.
\end{eqnarray*}
Applying the Markov property of $V$ at time $s_1$,   we have for $y\in
B$,
\begin{eqnarray*}
&&\P^y (V_{s_1+s_2}\in B(0, R-\tfrac12) \hbox{ and } \tau^V_B > s_1+s_2)\nonumber\\
&\geq & \P^y\left(\P^{V_{s_1}} \left( V_{s_2} \in B(0, R-\tfrac12 )
\hbox{ and } \tau^V_B>s_2 \right); \, V_{s_1}\in B(0, R-\tfrac14)
\hbox{
and } \tau^V_B>s_1 \right) \nonumber \\
&\geq & c_5 s_2 \P^y\left(  V_{s_1}\in B(0, R-\tfrac14 ) \hbox{ and
} \tau^V_B>s_1 \right) \nonumber\\
&\geq & c_5^2 s_2 s_1. \nonumber \\
\end{eqnarray*}
Repeating this at most $[4R]$ number of times, for any $y \in B$,
there exists $s_0 \le \tfrac14$ such that
\begin{equation}\label{P34a}
\P^y \left(V_{s_0}\in B(0, 1/4) \hbox{ and } \tau^V_B > s_0 \right)
\ge c_{11} (R-|y|)^\beta.
\end{equation}

Now let
$$ J_2(x,y) = J_0(x,y) + \kappa_2 |x-y|^{-d-\beta}\mathbbm{1}_{( |x-y|> 1)},$$
and $Y$ be the associated L\'evy process; note that $J_2$ satisfies
(\ref{ybeta}). We can use Remark \ref{meyer} to construct $Y$ from
$V$. Let $T$ be the time of the first added jump, so that $V_s=Y_s$,
$0\le s<T$ and $|\Delta Y_T|>1$. Let $B'=B(0,\half)$. If $x\in B'$
and $A \subset B'$ then $T \ge \tau^Y_{B'}\ge \tau^V_{B'}$.
Therefore
\[  \P^x \left( V_t\in A \hbox{ and } \tau_{B'}^V >t \right) \ge
\P^x \left( Y_t\in A \hbox{ and }  \tau_{B'}^Y >t \right),
 \]
  from which it follows that
$$ v^{B}(t,x,y) \ge  v^{B'}(t,x,y) \ge p_Y^{B'}(t,x,y), $$
for $x,y \in B'$.
 So, using Lemma \ref{stabproc}
\begin{equation} \label{vmidlb}
v^{B}(t,x,y) \ge c_{12}  \q \text { for } x, y \in B(0,\tfrac14)
\hbox{ and }  \tfrac34 \le t \le t_1.
\end{equation}
Hence for $t\in [t_0,t_1]$, $y\in B$ and $z\in B(0, \tfrac14)$, with
$s_0$ the time corresponding to $y$ in (\ref{P34a}), we have from
  (\ref{vmidlb})
\begin{align*}
 v^B(t,y,z)
 &\ge \int_{B(0, 1/4)} v^B(s_0,y,w)v^B(t-s_0,w,z) dw \\
 &\ge c_{12} \int_{B(0, 1/4)} v^B(s_0,y,w) dw\\
&= c_{12} \, \P^y \left( V_{s_0}\in B(0, 1/4) \hbox{ and }
\tau^V_B
>
 \tfrac14 \right)\\
 &\ge c_{12}  c_{11} (R-|y|)^\beta,
\end{align*}
which proves (\ref{vlowerb}). \qed

Define
\begin{equation}\label{eqn:a1}
 \vp(x)=\left( \left(R^2 -|x|^2 \right)^+ \right)^{12/(2-\beta)}.
\end{equation}
The following corollary follows immediately from Proposition \ref{utdecay}.

\begin{corollary}\label{vpovery}
For every $y\in B$, $\delta\in (0, 1)$ and
$\gamma \in \left(\frac{2-\beta}6, \, 1 \right]$,
there is a constant $C=(y, \beta, \delta,\xi)>0$ such that
\[
\vp(x)^\gamma/p^B(t,x, y)\leq C, \qquad  \hbox{for every }\, t\in [\delta,4]
\hbox{ and } \,  x \in B. \]
\end{corollary}

\bigskip
Recall  the definition of $\sF^B$ from Section 2.

\begin{lemma}\label{DFufinite}
For every $t>0$ and $y\in B$,  $p^B(t, x, y)$ as a function
of $x\in B$ is in $\sF^B$.
\end{lemma}

\proof  Fix $y\in B$. By Lemma \ref{delut}, if $t>0$,
$f(x):=p^B(t/2, x, y)\in L^2(B, dx)$ and so (see \cite[Lemma 1.3.3]{FOT94}),
$p^B(t, x, y)=P^B_{t/2} f(x)$ as a function of $x\in B$
is in $\sF^B$.
\qed

\begin{lemma}\label{DFvputfinite}
For each $t>0$ and $y_0\in B$,
the function $\vp (\cdot )/p^B(t, \cdot, y_0)$ is in $\sF^B$.
\end{lemma}

\proof Without loss of generality, we may assume that $t<4$.
By Corollary \ref{vpovery},
$\vp (\cdot )^{1/2}/p^B(t, \cdot, y_0)$ is bounded on $B$.
By extending the function $x\mapsto \vp (\cdot )/p^B(t, \cdot, y_0)$
to be zero on $B^c$, we see that $x\mapsto u(x):= \vp (\cdot )/p^B(t, \cdot, y_0)$
vanishes continuously on $B^c$.
Recall the killing function $\kappa_B$ from (\ref{eqn:c5}).
It is clear that $\int_Bu(x)^2 (dx)<\infty$.
So by (\ref{eqn:c3})--(\ref{eqn:c4}) and
 Lemma \ref{init-m-bounds},
 it suffices to show
\begin{equation} \label{vputa}
\int_{B\times B} \left( \frac{\vp (x )}{p^B(t, x, y_0)}
-\frac{\vp (y )}{p^B(t, y, y_0)} \right)^2 |x-y|^{-d-\beta} dx\, dy<\infty.
\end{equation}
We look at
\[
\int\int_{x,y\in B,\vp(y)\geq \vp(x)} \Big(\frac{\vp(y)}{p^B(t, y, y_0)}
-\frac{\vp(x)}{p^B(t,x, y_0)}\Big)^2 |x-y|^{-d-\beta} dx\, dy;
\]
if we show this is finite, we will have finiteness of the
the integral over $\vp(x)>\vp(y)$ by the same argument, and combining
the two estimates shows (\ref{vputa}).
We need to bound
\begin{align*}
\int&\int_{x,y\in B, \vp(y)\geq \vp(x)}
\Big(\frac{\vp(y)-\vp(x)}{p^B(t, y, y_0)}\Big)^2 |x-y|^{-d-\beta}dx\,dy\\
& ~~~~+ \int\int_{x,y\in B, \vp(y)\geq \vp(x)} \Big( \vp(x) \frac{p^B(t, y, y_0)-p^B(t, x, y_0)}
{p^B(t, x, y_0)p^B(t, y, y_0)}\Big)^2 |x-y|^{-d-\beta}dx\,dy\\
&~~~~=I_1+I_2.
\end{align*}
On the set $\{(x, y)\in B\times B: \, \vp (y)>\vp (x)\}$,
since  $\vp(x)\leq \vp(x)^{1/2} \vp (y)^{1/2}$ and
$\vp (x)^{1/2}/p^B(t, x, y_0)$ is bounded in $x$ by Corollary \ref{vpovery},
the finiteness of $I_2$  follows from Lemma \ref{DFufinite}.
To handle $I_1$, let $\gamma =(2-\beta)/3$.
Note that $|\vp(y)-\vp(x)|^{2\gamma}\leq |\vp(y)|^{2\gamma}$,
and that $\vp (x)^{2\gamma}/ p^B(t, x, y_0)^2$ is bounded in $x$ by Lemma \ref{utdecay}.
Since $\vp$ has a bounded derivative,
$|\vp(y)-\vp(x)|^{2-2\gamma}\leq c_1|y-x|^{2-2\gamma}$.
As $|y-x|^{2-2\gamma-d-\beta}$ is integrable over $B\times B$ since
$2\gamma<2-\beta$, the finiteness of $I_1$ follows.
\qed

\begin{lemma}\label{Gprime}
Fix $y_0\in B$ and
let $G(t)=\int_B \vp(x) \log p^B(t, x, y_0) \, dx.$
Then for every $t>0$,
\[
G'(t)=-\sE \left(p^B(t, \cdot , y_0), \,
\frac{\vp (\cdot )}{ p^B(t, \cdot, y_0)} \right).
\]
\end{lemma}

\proof Write $(f,g)$ for $\int_B f(x)g(x)\, dx$.
Using Lemmas \ref{DFufinite} and \ref{DFvputfinite}, we know
by \cite[Lemma 1.3.4]{FOT94} that
\begin{align*}
-\sE \left(p^B(t, \cdot , y_0), \,\frac{\vp (\cdot )}{ p^B(t, \cdot, y_0)}\right)
&=\lim_{h\downarrow 0} \frac{1}{h} \left(P^B_hp^B(t, \cdot , y_0)-p^B(t, \cdot , y_0), \,
\frac{\vp (\cdot )}{ p^B(t, \cdot, y_0)} \right)\\
&=\lim_{h\downarrow 0} \frac{1}{h} \left(p^B(t+h, \cdot , y_0), -p^B(t, \cdot , y_0), \,
\frac{\vp (\cdot )}{ p^B(t, \cdot, y_0)}\right)\\
&=\lim_{h\downarrow 0} \frac{1}{h}\int_B \vp (x)
\left( \frac{p^B(t+h, x , y_0)}{p^B(t, x, y_0)}-1 \right) dx.
\end{align*}
Also
$$
G'(t)=\lim_{h\to 0} \frac1{h}\int_B (\log p^B(t+h, x , y_0)-\log p^B(t, x , y_0)) \vp (x) dx.
$$
Let
$$ F(h)=\left[\log p^B(t+h, x , y_0)-\log p^B(t, x , y_0) -
\left( \frac{p^B(t+h, x , y_0)}{p^B(t, x, y_0)}-1 \right)
\right]\vp(x).
$$
Then
\begin{align*}
F'(h)&=\frac{\partial}{\partial t} p^B(t+h,x,y_0)\left(
\frac{1}{p^B(t+h,x,y_0)}-\frac{1}{p^B(t,x,y_0)}\right)\vp(x)\\
&=\frac{\partial}{\partial t} p^B(t+h,x,y_0)(p^B(t+h,x,y_0)-p^B(t,x,y_0)
\frac{\vp(x)}{p^B(t+h,x,y_0) p^B(t,x,y_0)}.
\end{align*}
By the mean value theorem,
$F(h)/h=F'(h^*)$ for some $h^*=h^*(x,y_0,h)\in (0,h)$.
Hence by Lemma \ref{delut} and Corollary \ref{vpovery},
$F(h)/h$ tends to 0 uniformly in $x\in B$ as $h\to 0$.
The lemma now follows from the dominated convergence theorem.
\qed

We need a weighted Poincar\'{e} inequality, which we derive along the lines of
the appendix to \cite{SCS}.

\begin{proposition}\label{wtedPI}
Let $R\in [1,4]$, $y_0\in \R^d$, $B=B(y_0,R)$,
\begin{equation}\label{defphi}
\phi_R(x)=c_1\left(R^2-|x-y_0|^2 \right)^{12/(2-\beta)} \,
\mathbbm{1}_{B}(x)
\end{equation}
with normalizing constant $c_1>0$ chosen so that $\int_B \phi_R (x) dx=1$,  and
set
\[ \overline f= \int_{B} f(x)\phi_R(x)\, dx . \]
There exists a constant $c_2$ depending on $R$ but not $f$ or $y_0$ such that
\[ \int_B |f(x)-\overline f|^2 \phi_R(x) dx
\leq c_2 \int_B \int_B (f(x)-f(y))^2 \phi_R(x)\land \phi_R(y)\, J(x,y) \, dx\, dy.
\]
\end{proposition}

\proof
If $B$ is any ball, let
\[ \sE_B(f,f)=\int_B\int_B (f(y)-f(x))^2 J(x,y)\, dx\, dy. \]
Set $f_B=|B|^{-1}\int_B f(x) dx$. If $B$ is any ball of radius $r\leq 1$, using (A3) we have
\begin{align*}
\int_B |f(x)-f_B|^2 dx&=\int_B (f(x)^2 -(f_B)^2) dx\\
&=\tfrac12|B|^{-1} \int_B\int_B (f(x)-f(y))^2 \, dx\, dy\\
&\leq \tfrac12 |B|^{-1}(\kappa_1 r^{-d-\al})^{-1}
  \int_B\int_B (f(x)-f(y))^2 J(x,y)\,dx\,dy\\
&= c_1 r^\al \sE_B(f,f).
\end{align*}
We now follow the proof in the appendix of \cite{SCS} closely, with the principal
changes being to use $\sE_B(f,f)$ in place of $\int_B |\nabla f(x)|^2 dx$,
$\int_B\int_B (f(x)-f(y))^2 (\phi(x)\land \phi(y)) \, J(x,y)\, dx\, dy$ in place
of $\int_B |\nabla f(x)|^2 \phi(x)\, dx$, and $r^\al$ in place of $r^2$.
\qed

\bigskip

\begin{proposition}\label{lowerb-prop}
Let $J$ satisfy the conditions (A1)--(A3). Let $\xi\in (0,1)$,
$J_\xi$ be defined by (\ref{Jxidefine}), and $X^{(\xi)}$ be the
Hunt process associated with $(\sE^{(\xi )}, \sF^{(\xi)})$.
Let  $y_0\in \R^d$ and  $\delta\in (0,1/2)$.
Let $R\in [1,4]$, $B=B(y_0,R)$, and
$p_\xi^B(t,x,y)$ be the transition density of $X^{(\xi)}$
killed on exiting $B$.
Then there exists
a positive constant $C$ that depends on
$\al, \beta, \kappa_1, \kappa_2$, and $\delta$, but not on $\xi$ or $y_0$
such that for all $t\in [\delta, 2]$
\begin{equation} \label{lowerbound-a}
p_\xi^B(t,x,y) \geq C  \end{equation}
for every $(x,y)\in
 B(y_0, 3R/4)\times B(y_0, 3R/4)$.
\end{proposition}

\proof
By a change of coordinate systems, without loss of generality
we may assume that
$y_0=0$ and so $B=B(0,R)$. Fix an arbitrary $x_0\in B(0,3R/4)$ and write
$u(t,x)=p_\xi^B(t,x_0,x)$.
Let  $\phi :\R^d \to \R_+$ be equal to  $\phi_R$ as defined in (\ref{defphi}).
Set for $t \in (0,\infty)$
\begin{align*}
 r(t,x)&=u(t,x)/\phi(x)^{1/2},\\
H(t) &= \int_B \phi(y) \log u(t,y) \; dy,\\
G(t) &= \int_B \phi(y) \log r(t,y) \; dy =\int_B \phi(y) \log u(t,y) \, dy-c_1=H(t)-c_1.
\end{align*}

Then by Lemma \ref{Gprime}
\begin{equation}  \label{G-Eeq}
G'(t)  =- \sE^{(\xi)} \Big( u(t, \cdot), \frac{\phi}{u(t, \cdot)} \Big) \;.
\end{equation}

The reason we work with $\sE^{(\xi)}$ rather than $\sE$ is so that we can use
 Lemma \ref{Gprime} to obtain (\ref{G-Eeq}). The remainder of the argument
does not use the condition (A4)($\xi$), and in particular the constants
can be taken to be independent of $\xi$.

By (\ref{eqn:c4}),
\begin{align*}
 G'(t)&=-\int_B\int_B \frac{[u(t,y)-u(t,x)]}{u(t,x)u(t,y)}
 [u(t,x)\phi(y)-\phi(x)u(t,y)] J^{(\xi)}(x,y) \,dy\,dx\\
&\qquad  -\int_B \vp (x) \kappa_B(x) dx.
\end{align*}
The main step is to show that for all $t$ in (0,1] one has
\begin{equation}\label{firststep}
G'(t) \geq c_2 \int_B\int_B [\log u(t,y)-\log u(t,x)]^2
(\phi(x)\land \phi(y)) J^{(\xi)}(x,y)\, dx\, dy -c_3
\end{equation}
for positive constants $c_2, c_3$.

\medskip

Setting $a=u(t,y)/u(t,x)$ and $b=\phi(y)/\phi(x)$, we see that
\begin{align}
\frac{[u(t,y)-u(t,x)]}{u(t,x)u(t,y)}& [u(t,x)\phi(y)-\phi(x)u(t,y)]\notag\\
&=\phi(x)\Big(b-\frac{b}{a}-a+1\Big)\label{Galg}\\
&= \phi(x)\Big[\Big((1-b^{1/2}\Big)^2-b^{1/2}\Big(
\frac{a}{b^{1/2}}+\frac{b^{1/2}}{a}-2\Big)\Big]. \notag
\end{align}
Using the inequality
\[ A+\frac{1}{A}-2\geq (\log A)^2, \qquad A>0, \]
with $A=a/\sqrt{b}$,
the right hand side of (\ref{Galg}) is bounded above by
\[ (\vp(x)^{1/2}-\vp(y)^{1/2})^2 -(\phi(x)\land \phi(y))\,
  (\log r(t,y)-\log r(t,x))^2.  \]
Substituting in the formula for $G'(t)$ and using Proposition \ref{wtedPI}
and (\ref{eqn:c5}),
\begin{align*} H'(t)=G'(t)
 &\geq -c_4+\int_B\int_B (\log r(t,y)-\log r(t,x))^2 \big(\phi(x)\land \phi(y)\big)
J^{(\xi)}(x,y)\, dx\, dy\\
&\geq -c_4 +c_5\int_B (\log r(t,y)-G(t))^2 \phi(y)\, dy\\
&\geq-c_6+c_7\int_B (\log u(t,y)-H(t))^2 \phi(y) \, dy\;.
\end{align*}
In the first inequality we used the fact
\[ \int_B\int_B (\vp(x)^{1/2}-\vp(y)^{1/2})^2 J^{(\xi)}(x,y) \, dx\, dy\leq c_8,  \]
which follows from (A3).
Recall the constant $t_0$ from Proposition \ref{exittime}. We may assume that
$\delta \leq t_0$.
By Proposition \ref{exittime},  for every $t\leq t_0$,
\begin{align*}
\int_{B(x_0, 1/4)} u(t,x)\, dx
&\geq \bP^{x_0}\left( \sup_{s\in [0, t_0]} |X_s-X_0| <1/4 \right)\label{small3dev}\\
&\geq 1-\bP^{x_0}\left(\sup_{s\in [0, t_0]} |X_s-x_0|\geq 1/4 \right)
\geq \tfrac34 .
\end{align*}

Choose $K$ such that $|B(x_0,1/4)| e^{-K}= \frac14$ and define
$$D_t=\{x\in B(x_0, 1/4): u(t,x)\geq e^{-K}\}.
$$
By Theorem \ref{rough-above-lem}, if $t\leq t_0$
\begin{align*}
\frac34 \leq \int _{B(x_0, 1/4)} u(t,x) dx
 &= \int_{D_t} u(t,x) \; dx +\int_{B(x_0,1/4)\setminus D_t} u(t,x) \; dx\\
&\leq c_9|D_t|t^{-d/\alpha} +|B(x_0, 1/4)| e^{-K}.
\end{align*}
Therefore
\[ |D_t|\geq \frac{t^{d/\alpha}}{2c_9}\geq c_{10} >0  \]
if $t\in [\delta/4, t_0]$.
Note that the positive constant $c_{10}$ can be
chosen to be independent of the $\xi$ in condition \a4xi.

Jensen's inequality tells us  that if $t\leq t_0$
$$ H(t)= \int_B\big(\log u(t,x)\big) \phi(x)\, dx
\leq  \log \int_B u(t,x)  \phi(x)\, dx
\leq \log \| \vp\|_\infty :=\bar H .  $$

On $D_t$, $\log u(t, x)\geq -K$ so there are only four possible cases:
\begin{description}
\item{(a)} If $\log u(t,x)>0$ and $H(t)\leq 0$,
then $(\log u(t,x)-H(t))^2\geq H(t)^2$.
\item{(b)} If $\log u(t,x)>0$ and $0<H(t)\leq \bar{H}$,
then
$$(\log u(t,x)-H(t))^2\geq 0\geq H(t)^2-\bar{H}^2.$$
\item{(c)} If $-K\leq \log u(t,x)\leq 0$ and $|H(t)|\geq 2K$, then
$(\log u(t,x)-H(t))^2\geq \frac14 H(t)^2$.
\item{(d)} If $-K\leq \log u(t,x)\leq 0$
and $|H(t)|<2K$, then
$$ (\log u(t,x)-H(t))^2\geq0\geq \frac14 H(t)^2 -K^2. $$
\end{description}
Thus we conclude
$$
\left( \log u(t, x) -H(t) \right)^2 \geq \frac14 H(t)^2 - (\bar H \vee K)^2
\qquad \hbox{on } D_t.
$$

Since $\phi$ is bounded below by $c_{11}>0$ on $B(x_0,1/4)$, then
\begin{align*}
& c_7 \int_B (\log u(t,x)-H(t))^2 \phi(x) dx -c_6
  \geq c_7 \int_{D_t}  (\log u(t,x)-H(t))^2 \phi(x) dx -c_6\\
& \quad \geq c_{12} |D_t| \Big(\frac14 H(t)^2-(\bar{H} \lor K)^2\Big) -c_6.
\end{align*}
We therefore have
\[ H'(t)\geq F H(t)^2-E, \qquad t\in [\delta/4,t_0]  \]
for some positive constants $E$ and $F$ that are independent of $\xi$.

\medskip

Now we do some calculus. Let $t_2\in [\delta/2,\, t_0 \wedge 2]$.
Let $Q:=\max( 16E, (16E/F)^{1/2})$.
Suppose $H(t_2)\leq -Q$. Since $H'(t)\geq -E$,
\begin{align}
H(t_2)-H(t)\geq -2E \qquad \hbox{for } t\in [\delta/4,t_2], \label{change-eps}
\end{align}
This implies $H(t)\leq -Q/2$. Since $FQ^2/4\geq 4E$,
$E<\frac{F}{2} H(t)^2$ and hence
\[ H'(t)\geq \frac{F}{2} H(t)^2. \]
Integrating $H'/H^2\geq F/2$ over $[\frac{\delta}4, \, t_2]$ yields
\[ \frac{1}{H(t_2)}-\frac{1}{H(\delta/4)}
  \leq -\frac{F}{2} (t_2-\delta/4)\leq -\frac{F\delta}{8}. \]
Since $H(\delta/4)\leq -Q/2<0$, we have $1/H(t_2)\leq -F\delta /16$, that is,
\[ H(t_2)\geq -\frac{16}{F\delta}.  \]
This proves that either $H(t_2)\geq -Q$ or $H(t_2)\geq -16/(F\delta)$.
Thus in either case, $H(t_2)\geq -U$, where
$U=U(\delta):=\max\{Q, \,16/(F\delta)\}>0$, and so $G(t_2)\geq -U + c_1$.

\medskip
Now for every $x_0, x_1\in B(0, 3R/4)$, applying
the above first with $x_0$ and then with $x_0$ replaced by
$x_1$, we have
\begin{align*}\log p_\xi^B(2t_2,x_0,x_1)
  &=\log \int p_\xi^B(t_2,x_0,z)p_\xi^B(t_2,x_1,z)\, dz \\
&\geq \log
 \int_B p_\xi^B(t_2,x_0,z)p_\xi^B(t_2,x_1,z) \phi(z) \, dz -\log \| \phi \|_\infty\\
&\geq \int_B \log \Big(p_\xi^B(t_2,x_0,z)p_\xi^B(t_2,x_1,z)\Big) \phi(z)\, dz
   -\log \| \phi \|_\infty \\
&= \int_B \log p_\xi^B(t_2,x_0,z) \phi(z) dz +\int_B\log p_\xi^B(t_2,x_1,z) \phi(z) dz \\
  &~~~~~~~ -\log \| \phi\|_\infty \\
& \geq -2(U+c_{12}),
\end{align*}
that is, $p_\xi^B(2t_2,x_0,x_1)\geq e^{-2(U+c_{12})}$.
A repeated use of the semigroup property (but at most $2/t_2$ more times) then shows
$p_\xi^B(t,x_0,x_1)\geq c_{13}(\delta)$ for every
$t\in [\delta/2, \, 2]$.
\qed

\begin{theorem}\label{lowerb-theo-a}
Let $J$ satisfy the conditions (A1)--(A3). Let $\xi\in (0,1)$,
$J_\xi$ be defined by (\ref{Jxidefine}), and $X^{(\xi)}$ be the
Hunt process associated with $(\sE^{(\xi )}, \sF^{(\xi)})$.
Let  $y_0\in \R^d$ and  $\delta\in (0,1/2)$.
Let $R>0, T>1/2$, $B=B(y_0,R)$, and
$p_\xi^B(t,x,y)$ be the transition density of $X^{(\xi)}$
killed on exiting $B$.
Then there exists
a positive constant $C=C(R)$ that depends on
$\al, \beta, \kappa_1, \kappa_2$, $R$ and $\delta$, but not on $\xi$ or $y_0$
such that for all $t\in [\delta, T]$
\begin{equation} \label{lowerbound-a2}
p_\xi^B(t,x,y) \geq C  \end{equation}
for every
$(x,y)\in (B(y_0, 3R/4)\times (B(y_0, 3R/4)$.
\end{theorem}

\proof By a change of coordinate system, we assume without
loss of generality that $y_0=0$. Suppose first that $T \geq 2$ and $R\in [1,4]$.
Let $r=3R/4, n=[T]$, and $\delta=1/n$. So if $t\in [1,T]$, then
$t/n\in [\delta,2]$. By the semigroup property, if $x,y\in B(0,3R/4)\setminus \sN$, then
\begin{align}
p^B(t,x,y)
&\geq \int_{B(0,r)}\cdots \int_{B(0,r)} p^B(t/n,x,z_1)p^B(t/n,z_1,z_2)\cdots
\label{chaining}\\
&~~~~~~~~~~~~~~~~~~~~~~~~~~~~~ p^B(t/n,z_{n-2},z_{n-1})p^B(t/n,z_{n-1},y)\, dz_1\cdots
dz_{n-1}\nonumber\\
&\geq c_1|B(0,r)|^{n-1}\geq c_2 \nonumber
\end{align}
by Proposition \ref{lowerb-prop}.
We therefore have the conclusion of the theorem for all $T\geq 2$.

Next suppose $R>4$. Suppose $x,y\in B(0,3R/4)\setminus \sN$, $t\in [\delta, T]$,
$n=2[|x-y|]+1$, and let $z_0=x, z_1, \ldots, z_{n-1}, z_n=y$ be equally spaced
points on the line segment joining $x$ and $y$.
Then $|z_{i+1}-z_i|\leq 1/2\leq R/8$.
Set $r=1$. Using
(\ref{chaining}) and Proposition \ref{lowerb-prop} with $\delta$ replaced
by $\delta/n$, we again obtain our conclusion.

Finally, suppose $R<1$. Fix $\delta$ and $T$. Consider the process
$Z=R^{-1}X^{(\xi)}$ with corresponding jump kernel $J_Z$. By a change of
variables, we see that it suffices to obtain a lower bound on $p_Z^B(t,x,y)$
for $x,y\in B(0,3/4)$,
and $p_Z^B$ is the transition
density for $Z$ killed on exiting $B(0,1)$. The jump kernel corresponding
to $Z$ is $J_Z(x,y)=R^{-d} J(Rx,Ry)$. Let $J_{Z^{(0)}}(x,y)
=J_Z(x,y) \mathbbm{1}_{(|x-y|<1)}$. It is easy to see that $J_{Z^{(0)}}$ satisfies
(A1)--(A3) and \a4xi (but with different constants $\kappa_1, \kappa_2, \xi$).
Let $Z^{(0)}$ be the process corresponding to $J_{Z^{(0)}}$ and construct
$Z$ from $Z^{(0)}$ using Remark \ref{meyer}. Then if $A\subset B(0,3/4)$,
$\tau_B=\inf\{t: Z_t\notin B\}$, and $\tau^0_B=\inf\{t: Z^{(0)}_t\notin B\}$,
by independence, and using Lemma \ref{mey-bound},
\begin{align*}
\int_A p^B_Z(t,x,y)\, dy&=\P^x(Z_t\in A, \tau_B>t)\\
 &\ge e^{-(\sup \sJ)T} \int_A p^B_{Z^{(0)}}(t,x,y)\, dy
\geq c_3 |A|,
\end{align*}
which proves the theorem in this case as well.
\qed

\section{Mosco convergence}

In this section we first prove some general results on Mosco convergence.
We will then use these to prove Theorems \ref{lowerbound-theo} and
\ref{moscoconv}.

Let us first recall the definition of Mosco convergence and its properties.
Let $E$ be a locally compact separable metric space and $m$
a Radon measure on $E$ with full support.
Given a  densely defined quadratic form $(\sE, \sF)$ in $L^2(E; m)$,
we can extend its domain of definition to $L^2(E; m)$
by setting $\sE(u, u)=\infty$ for $u\in L^2(E; m)\setminus \sF$.
Throughout this section we will use this extension and, unless otherwise
specified, all the quadratic forms encountered will be assumed to be densely defined
in $L^2(E; m)$.
Recall that given $(\sE, \sF)$ we set $\sE_1(f,f)=\sE(f,f)+||f||_2^2$.

\begin{definition}\label{D:2.1}
A sequence of closed quadratic forms
$\{(\sE^{n}, \sF^{n})\}$
on $L^2(E; m)$ is said to be convergent
to a closed quadratic form $(\sE, \sF)$ on $L^2(E; m) $
in the sense of Mosco (cf. \cite{Mosco}) if
\begin{description}
\item{\rm (a)} For every sequence $\{u_n, n\geq 1 \}$ in $L^2(E; m)$
that converges weakly to $u$ in $L^2(E; m)$,
$$\liminf_{n\to \infty} \sE^{n}(u_n, u_n) \geq \sE (u, u) , $$
\item{\rm (b)} For every $u\in L^2(E; m)$, there is a sequence
$\{u_n, n\geq 1\}$ in $L^2(E; m)$ converging strongly to $u$
in $L^2(E; m)$ such that
$$ \limsup_{n\to \infty} \sE^{n}(u_n, u_n) \leq \sE (u, u) . $$
\end{description}
\end{definition}

Let $\{P_t, \, t\geq 0\}$ and $\{P^{n}_t, \, t\geq 0\}$
be the semigroups of  $(\sE, \sF)$ and $(\sE^{n}, \sF^{n})$, respectively,
and $\{G_\alpha, \alpha >0 \}$ and $\{G^{n} _\alpha, \alpha >0 \}$
their corresponding resolvents, respectively.
The following result is known (see Theorem 2.4.1 and
Corollary 2.6.1 of \cite{Mosco}).

\begin{proposition} \label{P:2.1}
Let $(\sE, \sF)$ and $\{(\sE^{n}, \sF^{n}), \, n\geq 1\}$
be closed quadratic forms on $L^2(E; m)$.
The following are equivalent:
\newline (a) $(\sE^{n}, \sF^{n})$ converges to $(\sE, \sF)$
in the sense of Mosco;
\newline (b)
For every $\alpha >0$ and $f\in L^2(E; m)$,
$G^{n}_\alpha f$ converges to
 $G_\alpha f$ in $L^2(E; m)$;
  \newline (c)
For every $t >0$ and $f\in L^2(E; m)$,
$P^{n}_t f$ converges to $P_t f$ in $L^2(E; m)$.
\end{proposition}

Here is a criterion for Mosco convergence to hold.

\begin{theorem}\label{T4.1}
Suppose
\begin{description}
\item{(i )} $\sF^{n}\subset \sF$ for every $n\geq 1$ and $\sE^{n}(u, u)
\geq \sE (u, u)$ for every $u\in \sF^{n}$.

\item{(ii)} There is a common  core ${\cal C}$ for the Dirichlet forms
$(\sE^{n}, \sF^{n})$ and $(\sE, \sF)$ such that
$$ \lim_{n\to \infty} \sE^{n}(u, u) = \sE (u, u) \qquad
\hbox{for every } u\in {\cal C}.
$$
\end{description}
Then $(\sE^{n}, \sF^{n})$ converges to $(\sE, \sF)$ in the
sense of Mosco.
\end{theorem}

\pf
Let $\{u_k, k\geq 1 \}$ be a sequence in $L^2(E; m)$
that converges weakly to $u$ in $L^2(E; m)$.
Without loss of generality,
we may assume that $\lim_{k\to \infty} \sE^{k} (u_k, u_k)$ exists and is finite.
This in particular implies that $u_k\in \sF^{k}\subset \sF$ for every $k\geq 1$ and
$\sup_{k\geq 1} \sE(u_k, u_k)<\infty$. Since $\{u_k, k\geq 1 \}$ is bounded in
$L^2(E; m)$, taking a subsequence if necessary, we may assume that
the Cesaro mean of $\{u_k, k\geq 1 \}$ converges in $(\sE_1, \sF)$ to
some function $v$ (see page 14 of \cite{Silver}).
As $u_k$  converges weakly to $u$, we must have $u=v\in \sF$.
Therefore,
$$
   \liminf_{k\to \infty} \sE^{k} (u_k, u_k)
   \geq \limsup_{k\to \infty }
\sE(u_k, u_k) \geq \limsup_{k\to \infty} \sE \left( \frac1k \sum_{j=1}^k u_j,
\frac1k \sum_{j=1}^k u_j \right)\geq \sE(u, u).
$$
The second inequality follows since the triangle inequality tells us
that
$$\sE \left(\frac{1}{k}\sum_{j-1}^k u_j, \frac{1}{k}\sum_{j=1}^k u_j
\right)^{1/2}\leq
\frac{1}{k} \sum_{j=1}^k \sE(u_j,u_j)^{1/2}.$$
This shows that the condition (a) of Definition \ref{D:2.1} is satisfied.

For any $u\in {\cal F}$, there
exists a sequence $\{v_j\}\subset \sC$ converging
strongly to $u$ in $L^2(E; m)$ such that
$$
\lim_{j\to \infty}{\cal E}(v_j, v_j)=
{\cal E}(u, u).
$$
Since for each $j\geq 1$,
$$ \lim_{n\to \infty} \sE^{n}(v_j, v_j) = \sE (v_j, v_j),
$$
using induction we can find an increasing subsequence
$\{n_j\}$ such that
$$
|{\cal E}^{n}(v_j, v_j)-{\cal E}(v_j, v_j)|\le 2^{-j}
\qquad \hbox{for } n\ge n_j.
$$
Put $u_1=\dots=u_{n_1-1}=0$ and
$u_{n_j}=\dots=u_{n_{j+1}-1}=v_j$ for $j\ge 1$.
It is easy to see that $\{u_k, k\geq 1\}$ is a sequence
in $\sC$ converging
strongly to $u$ in $L^2(E; m)$ such that
$$ \lim_{k\to\infty}{\cal E}^{k}(u_k, u_k)=
{\cal E}(u, u). $$
For $u\in L^2(E, m)\setminus \sF$, since
$\sE (u, u)=\infty$, it trivially holds that
$$ \limsup_{k\to \infty} \sE^{k} (u_k, u_k) \leq \sE(u, u). $$
This shows that the condition (b) in Definition \ref{D:2.1} is
satisfied. Hence we have shown that $(\sE^{k}, \sF^{k})$ is
Mosco-convergent to $(\sE, \sF)$. \qed

\bigskip

Let $X^k$ be the Hunt process associated with $(\sE^{k}, \sF^{k})$
and $X^{k, B}$ be the subprocess of $X^k$ killed upon exiting  an
open set $B$. It is known (see \cite{FOT94}) that the Dirichlet form
$(\sE^{k}, \sF^{k, B})$ of $X^k$ is given by
$$ \sF^{k, B}=\left\{ u\in \sF^{k}: \, u = 0 \ \ \sE^k\hbox{-q.e. on } B^c \right\}. $$

\medskip

\begin{theorem}\label{R:4.6}
Suppose $B$ is an open set and the following hold.
\begin{description}
\item{(i)} $\sF^k\subset \sF$ and $\sE^k(u, u)
\geq \sE (u, u)$ for every $u\in \sF^k$ and  every $k\geq 1$.

\item{(ii)} There is a common  core ${\cal C}$ for the Dirichlet forms
$(\sE^k, \sF^k)$ and $(\sE, \sF)$ such that
$$
\lim_{k\to \infty} \sE^k (u, u) = \sE (u, u) \qquad
 \hbox{for
every } u\in {\cal C}. $$ Furthermore, there is a common core ${\cal
C}_B\subset \sC$ for the Dirichlet forms $(\sE^k, \sF^{k, B})$ and
$(\sE, \sF^B)$.
\end{description}
Then $(\sE^k, \sF^{k,B})$ converges to $(\sE, \sF^B)$ in the sense
of Mosco.
\end{theorem}

\proof
To emphasize the domain of definition, for this proof only,
we write $\sE^B$ and $\sE^{k, B}$ for $(\sE^k, \sF^{k, B})$
and $(\sE, \sF^{B})$, respectively. With this notation,
$\sE^B(u, u)=\infty$ when $u\notin \sF^B$ and
$\sE^{k, B}(u, u)=\infty$ when $u\notin \sF^{k,B}$.

First note that, since $\sE^k_1 (u, u)\geq \sE_1 (u, u)$,
$$ \sF^k\subset \sF \qquad \hbox{and} \qquad \sF^{k, B} \subset \sF^B. $$
For any $v_k$ that converges weakly to $v$ in $L^2(B; dx)$,
we claim that
\begin{equation}\label{eqn:c6}
\liminf_{k\to\infty}{\cal E}^{k, B}(v_k, v_k)\ge
{\cal E}(v, v).
\end{equation}
Suppose that the left hand side of (\ref{eqn:c6}) is finite.
Then there is a
subsequence $\{n_k\}$ such that
$$\lim_{k\to \infty} \sE^{n_k, B}(v_{n_k}, v_{n_k})
= \liminf_{k\to\infty}{\cal E}^{k, B}(v_k, v_k) \quad \hbox{ and }
\quad \sup_{k\geq 1} \sE^{n_k, B} (v_{n_k}, v_{n_k}) < \infty.
$$
In particular, this implies that $v_{n_k} \in \sF^{n_k, B}$ and
$$ \sup_{k\geq 1} \sE^{B}_1 (v_{n_k}, v_{n_k}) \leq
\sup_{k\geq 1} \sE^{n_k, B}_1 (v_{n_k}, v_{n_k}) < \infty.
$$
By taking a subsequence if necessary, the Cesaro mean of
$\{v_{n_k}, k\geq 1\}$ converges in
in $\sF^B$ with respect to the Hilbert norm $\sqrt{\sE^B_1}$ to
a function $w$, which has to be $v$. This implies that
$v\in \sF^B$.
By extending $v_{n_k}$ and $v$ to take the value zero off $B$,
we have $v_{n_k} \in \sF^{n_k}$ and $v\in \sF$.
By Theorem \ref{T4.1},
$({\cal E}^{k}, {\cal F}^{k})$
converges to $({\cal E}, {\cal F})$ in the sense
of Mosco, and we have in particular that
$$
\liminf_{k\to\infty}{\cal E}^{k}(v_k, v_k)\ge
{\cal E}(v, v).
$$
As $\sE^{k, B}$ and $\sE^B$ agrees with $\sE^k$ and $\sE$ on $\sF^{k, B}$ and
$\sF^B$, respectively, this proves (\ref{eqn:c6}).

Noting that $\sC_B$ is a common core
for $({\cal E}, \sF^B)$ and $({\cal E}^k, \sF^{k, B})$,
it can be shown that the second condition (b) in Definition \ref{D:2.1}
holds for $({\cal E}^{k, B}, {\cal F}^{k, B})$ and $({\cal E}^B, {\cal F}^{B})$
in much the same way as that in the proof of Theorem \ref{T4.1}.
\qed

\bigskip

\noindent{\bf Proof of Theorem \ref{moscoconv}:}
Let $\delta_k$ be a sequence of positive numbers decreasing to 0.
Set
\begin{equation}\label{Jkdefine}  J_k(x, y)=\
\begin{cases}
J(x, y) &\hbox{for } |x-y|\geq \delta_k; \\
 \kappa_2 |y-x|^{-d-\beta} &\hbox{for }  |x-y|<\delta_k,
\end{cases}
\end{equation}
and define $(\sE^k, \sF^k)$ in the same way as we defined  $(\sE, \sF)$
in  (\ref{Dirform})--(\ref{eqn:domain}).
Note that $\sE^k$ satisfies A4($\delta_k$).
Take $\sE^k(f,f)=+\infty$ if $f \in L^2(\bR^d,dx)\setminus \sF^k$.

It is clear that $J_k(x, y)$ decreases to $J(x, y)$ as $k\uparrow
\infty$, and so $\sF^k\subset \sF$ and $\sE^k(u, u)\geq \sE(u, u)$
on $\sF^k$ for every $k\geq 1$. By Lemma \ref{init-m-bounds},
$$ \sF^k =  \H^{\beta/2}(\R^d). $$
and thus in particular $\sF^k$ is independent of $k$.
Note that $C^1_1(\R^d)$ is the common core of $( \sE^k, \sF^k)$ for
$k\geq 1$ and for $(\sE, \sF)$, and that
 $$  \lim_{k\to \infty} \sE^k (u, u) = \sE(u, u) \qquad
 \hbox{for every } u\in C^1_c(\R^d).
 $$
  Theorem
\ref{moscoconv} now follows from Theorem \ref{T4.1}. \qed

\begin{corollary}\label{C4.2}
Let $B$ be a ball and define $\sF^{k,B}$ by (\ref{eqn:c3}), where
$(\sE^k, \sF^k)$ is as in the proof of Theorem \ref{moscoconv}.
Then $(\sE^k, \sF^{k,B})$ converges in the sense of Mosco to
$(\sE, \sF^B)$.
\end{corollary}

\proof
Note that $C^1_c(B)\subset C^1_c(\R^d)$ is a common core for $(\sE,
\sF^B)$ and $(\sE^k, \sF^{k, B})$.
 The conclusion of the corollary follows directly from Theorem \ref{R:4.6} and
 the proof of  Theorem \ref{moscoconv}. \qed

\medskip

\noindent {\bf Proof of Theorem \ref{lowerbound-theo}}:
Let $\delta_k$ be a sequence of positive numbers decreasing to 0.
Define $J_k(x,y)$ by (\ref{Jkdefine})
and define $(\sE^k, \sF^k)$ as in the proof of Theorem \ref{moscoconv}.
Clearly $J_k$ satisfies the conditions (A1)--(A3) as well as
(A4)($\delta_k)$
with the same $\kappa_1$ and $\kappa_2$ as $J$.

Let $p^{k, B}(t, x, y)$ and $p^B(t, x, y)$ denote
the transition density functions of $X^{k, B}$ and $X^B$
respectively. It follows from Theorem \ref{lowerb-theo-a} and
Proposition \ref{P:2.1} that for any given $\delta \in (0, 1)$,
there is a constant $c=c(\delta) >0$ such that
for any bounded non-negative functions
$f$ and $g$ on $B$ and $t\in [\delta, 2]$,
\begin{align*}
 \int_{B(y_0, 3R/4)\times B(y_0, 3R/4)} & p^B(t, x, y) f(x)g(y) dx\,dy \\
&=\lim_{k\to \infty }\int_{B(y_0, 3R/4)
   \times B(y_0, 3R/4)} p^{k,B}(t, x, y) f(x)g(y) dx\,dy \\
&\geq c \int_{B(y_0, 3R/4)\times B(y_0, 3R/4)} f(x)g(y) dx\,dy .
\end{align*}
This implies that $p(t, x, y)\geq c$ for almost every $ x, y \in
B(y_0, 3R/4)$. On the other hand, it follows from the proof of
Theorem \ref{rough-above-lem} that there is a properly exceptional
set $\sN$ so that $p(t, x, y)$ is well-defined on
$(B\setminus \sN )\times (B\setminus \sN)$ and that for each
fixed $y\in B\setminus \sN$,
$x\mapsto p^B(t, x, y)$ is $X^B$-quasi-continuous (and hence
$X$-quasi-continuous).
It follows
 that $p(t, x, y)\geq c$ for every $t\in [\delta, 2]$ and
every $x, y \in B (y_0, 3R/4)\setminus \sN$. \qed

\section{Parabolic Harnack inequality}

In this section we prove Theorem \ref{Maintheorem} in the case
$R\geq 1$. The argument uses balayage; see \cite{blumgetoor},
Chapter VI, for details.

\medskip
\noindent {\bf Proof of Theorem \ref{Maintheorem}:}
Without loss of generality we may assume the following: by a change of
coordinate system, we may assume $x_0=0$; by the Markov property we
may assume $t_0=0$; by looking at the process $X'_t=X_{t/T}$, we see
that the jump kernel corresponding
to $X'$ satisfies (A1)--(A3) (but with different $\kappa_1, \kappa_2$),
so we may assume $T=1$.
With these assumptions $Q=(0,5)\times B(0,4R)$. Recall the notation for
hitting and exit times given in (\ref{hittingtimedef}).
Let
$E=(\frac12,\frac92)\times B(0,3R)$, $D=(\frac34,\frac{17}{4})\times B(0,2R)$, and write
$A=\overline{B(0,3R)}\setminus B(0,2R)$ and  $B=B(0,4R)$.
 By the martingale
property,
\begin{equation}
P^B_{t-s} u(s,x)\leq u(t,x), \q \hbox { for  $s<t$ with $(s,x)$, $(y, t) \in Q$}.
\end{equation}
This says that the function $u$ is excessive with respect to the
space-time subprocess $(V^Q, X^Q)$ of $(V, X)$ killed upon exiting
$Q$, where
$V_t=V_0-t$. We can define $u_E$, the r\'eduite of $u$ with respect
to $E$, by
$$ u_E(s,x)=\E^{(s,x)} [u(V_{T_E},X_{T_E}); T_E<\tau_Q].$$
The function $u_E$ is again excessive with respect to
the killed process $(V^Q, X^Q)$, is 0 on $Q^c$,
and agrees with $u$ on $E$; see \cite{blumgetoor}. Note
that the process $(t,X_t)$ is in duality with the process $(V_t,
X_t)$ in the sense of \cite[Chapter VI]{blumgetoor}. By the Riesz
decomposition theorem (cf. \cite[Theorem VI.2.11]{blumgetoor}),
$u_E$ is the potential of a measure $\nu_E$ supported on $\ol E$. This
means that if $(t,x) \in Q$ then
\begin{equation}\label{u-rep-a}
 u_E(t,x) =\int_E p^B(t-r,x,z)\, \nu_E(dr,dz).
\end{equation}
Here we have $p^B(s,x,y)=0$ if $s<0$.

Since the jumps of the process $X$ are bounded by 1,
$u_E$ is caloric on $(1/2, 9/2 )\times B(0, 2R)$.
It follows that the support of $\nu_E$ is contained in
$\ol E\setminus \left( (1/2, 9/2 )\times B(0, 2R) \right)$.
For $t\in (1/2, 5)$, let
\begin{align*}
  F_1(t) &=  [1/2,t]\times (A\setminus \sN),  \\
  F_2(t) &= \{\tfrac12\} \times \big(\ol {B(0,2R)} \setminus \sN \big), \\
  F(t) &= F_1(t) \cup F_2(t).
\end{align*}
Thus if $(t,x) \in D$ then we can write (\ref{u-rep-a}) as
\begin{equation}\label{u-rep-b}
u(t,x) = u_E(t,x) = \int_{F(t)}  p^B(t-r,x,z)\, \nu_E(dr,dz).
\end{equation}
Since $\nu_E$ is an equilibrium measure (i.e., a capacitary
measure), it does not charge polar sets; in particular, it does
not charge $[0,5]\times \sN$.

Consider (\ref{u-rep-b}) when
$(t_1,x_1) \in Q_{-}=(1, 2)\times (B(0,R)\setminus \sN)$.
If $(r,z)\in F_1(t)$ then
$|x_1-z|\ge R$, and thus by Proposition \ref{exittime}(a),
\begin{equation}\label{u-rep-c}
  p^B(t_1-r,x_1,z) \le  p(t_1-r,x_1,z) \le c_1.
\end{equation}
If $(r,z)\in F_2(t)$, then
$t_1-r\ge \frac12$ and by Theorem \ref{rough-above-lem} again
(\ref{u-rep-c}) holds.

Now let $(t_2,x_2) \in Q_{+}=(3, 4)\times (B(0, R)\setminus \sN)$.
If $(r,z) \in F(t_1)$ then $t_2-r\ge t_2-t_1 \ge 1$, and
$|x_2-z| \le 4$, so by
Theorem \ref{lowerbound-theo}
\begin{equation}\label{u-rep-d}
 p^B(t_2-r, x_2, z) \geq c_2.
\end{equation}
Hence
\begin{align*}
u(t_2, x_2) &=  \int_{F(t_2)} p^B(t_2-r,x_2,z) \nu_E(dr,dz)
\\
&\geq \int_{F(t_1)} p^B(t_2-r,x_2,z) \nu_E(dr,dz) \\
&\geq c_2 \nu_E(  F(t_1)) \\
&\geq (c_2/c_1) \int_{F(t_1)} p^B(t_1-r,x,z) \, \nu_E(dr,dz)\\
&=  (c_2/c_1)\,  u(t_1, x_1),
\end{align*}
giving the parabolic Harnack inequality with constant
$C= c_1/c_2$. \qed

\section{Harmonic functions need not be continuous}

One of the applications of scale invariant Harnack inequalities is
that they imply regularity, e.g., H\"older continuity of harmonic and
caloric functions, and resolvents. This can be used in order to remove
properly exceptional sets. It is therefore interesting to see that
such regularity can fail, even when a Harnack inequality holds.  We
say a function $h$ is harmonic in a domain $D\subset \R^d$ if
$h(X_{t\land \tau_{D_1}})$ is a right continuous martingale for every
subdomain $D_1$ with $\overline D_1\subset D$, where
$\tau_{D_1}=\inf\{t>0: X_t\notin D_1\}$.

In this section we construct a class of symmetric jump processes
satisfying our hypotheses where there exist bounded harmonic
functions that are not continuous. An interesting side result
related to this example is that the martingale problem for variable
order jump processes is not always well posed. See Remark
\ref{rem2}(d) and the results in \cite{HuKa05} for other examples
which are similar but lead to continuous harmonic functions and
Feller semigroups.

\medskip
Using the integral conditions given
in Theorems 11.2 and 11.5 of \cite{fris} we obtain:

\begin{lemma}\label{stablelim}
Let $X_t$ be a one-dimensional stable process of index $\alpha \in (0,1)$
and $\eps>0$. Then
\[  \liminf_{t\to 0} \frac{|X_t|}{t^{(1/\al)+\eps}} =\infty
 \qquad \hbox{and} \qquad
\limsup_{t\to 0} \frac{|X_t|}{t^{(1/\al)-\eps}} =0, \qquad \bP^0 \mbox{-a.s.} \]
\end{lemma}

 Before constructing the main counterexample on $\R^2$,
 we need to look at an auxiliary process $Y$. Let $0<a<b<2$ and set
for $z \in \R^2$, $z_1\ne z_2$ :
\begin{align}
 m(z_1,z_2) =
\begin{cases}
    \min(|z_1|^{-a-2}, |z_2|^{-b-2})   &\hbox{ if } |z_1|\vee |z_2| \le 1, \\
    0                         &\hbox{ if } |z_1|\vee |z_2| > 1. \\
\end{cases} \label{eq:m-choice}
\end{align}
Assume $|z_1|\vee |z_2| \le 1$. Note that
\begin{align*}
\frac{1}{|z_1|^{a+2}+|z_2|^{b+2}}
\le \min\Big(\frac{1}{|z_1|^{a+2}}, \frac{1}{|z_2|^{b+2}}\Big)
\leq \frac{2}{|z_1|^{a+2}+|z_2|^{b+2}} \,.
\end{align*}
This implies
\begin{align}
c_0 |z|^{-a-2} \leq \min(|z_1|^{-a-2}, |z_2|^{-b-2}) \leq c_1 |z|^{-b-2} \;,
\end{align}
where $c_0, c_1$ are independent of $z$. Now for $x, y \in \R^2$, $x \ne y$ let
\begin{align*}
J_0(x,y)=
\begin{cases}
m(|x_1-y_1|,|x_2-y_2|) \quad & \text{ if } x-y \in [-1,1]^2 \,, \\
0 & \text{ if } x-y \notin [-1,1]^2 \,.
\end{cases}
\end{align*}
Furthermore, set
\begin{equation} \label{alpha-beta}
\al = \frac{(a+1)(b+1)-1}{b+2}, \q \beta = \frac{(a+1)(b+1)-1}{a+2} \;.
\end{equation}

The following facts are now obvious.

\begin{lemma} \label{indices}
\begin{description}
\item{(a)} $J_0(x, y)$  is symmetric in $(x, y)$ and
\begin{equation}\label{a3ce}
 c_1 |x-y|^{-a-2} \le J_0(x,y) \le c_2  |x-y|^{-b-2}  \qquad
 \hbox{for } |x_1-y_1|\vee |x_2-y_2| \leq 1.
 \end{equation}

\item{(b)} $\al< \beta$.

\item{(c)} $a < (2-b)/b$  if and only if  $\beta
<1$.
\end{description}
\end{lemma}

Now choose $a$, $b$ with  $0<a<b<2$ and  $a < (2-b)/b$, so that
$0<\al<\beta<1$. Let $Y_t=(Y_t^1,Y_t^2)\in \bR^2$ be the pure jump
symmetric L\'evy process with jump intensity kernel $J_0$. The
following lemma explains the behavior of the marginals
$Y_t^1,Y_t^2$.

\begin{lemma} \label{easyint} We have
$$ n_1(z_1) :=\il_{-1}^1 m(z_1,z_2) dz_2 = \frac{2}{b+1} +  \frac{2b}{b+1}
 |z_1|^{-\al-1}
  \qq \hbox{for }  z_1 \in [-1,1], $$
$$ n_2(z_2):=\il_{-1}^1 m(z_1,z_2) dz_1 = \frac{2}{a+1}
  +  \frac{2a}{a+1} |z_2|^{-\beta-1}  \qq \hbox{for } z_2 \in [-1,1]. $$
\end{lemma}

\begin{proof}
\begin{align*}
\frac12  \il_{-1}^1 m(z_1,z_2) dz_2
&= \il_{0}^{|z_1|^\frac{a+2}{b+2}} |z_1|^{-a-2} dz_2
     + \il_{|z_1|^\frac{a+2}{b+2}}^1   |z_2|^{-b-2} dz_2 \\
&= |z_1|^{-\frac{(a+2)(b+1)}{b+2}} + \frac{1}{b+1}
 - \frac{1}{b+1} |z_1|^{-\frac{(a+2)(b+1)}{b+2}} \,.
\end{align*}
Note that $\frac{(a+2)(b+1)}{b+2} = \alpha +1$. The first assertion of
the lemma follows. The second one is proved analogously. \qed
\end{proof}

The coordinate processes $Y^i$ are one-dimensional L\'evy processes
with jump measure $n_i$; note however that $Y^1$ and $Y^2$ are
not independent.

Let
$$  V(\lam)=\{(x_1,x_2): |x_1|<\lam |x_2|\}, \q
  \tau_{V(\lam)} = \inf\{t >0: Y_t \not\in V(\lam)\}. $$
and write $V=V(1)$.

\begin{lemma} \label{inW} Let $\lam>0$;
then $\bP^0( \tau_{V(\lam)}>0)=1$.
\end{lemma}

\proof By Lemma \ref{easyint}, $n_1$ differs from the jump measure
of a stable process of index $\al$ by a finite measure. Therefore
$Y^1$ has the same local behavior at time 0 as a stable process
of index $\al$;
similarly $Y^2$ has the same local behavior at time 0 as a
stable process of index  $\beta$. (Note that points are polar for these two
processes.) Choose $\eps$ such that $\eps+ \tfrac1{\beta} <
\tfrac1{\al} -\eps$. By Lemma \ref{stablelim} for all sufficiently
small $t>0$,
\[ 0< |Y_t^1|\le t^{(1/\alpha)-\eps}\le\half \lam  t^{(1/\beta)+\eps}
\leq \half \lam |Y_t^2|. \]
\qed

\begin{remark}\label{R7.1}
{\rm
By Lemma \ref{inW} we can find $t_0>0$ such that
$\bP^0( \tau_{V(1/3)} \leq t_0 ) <  1/20$.
 Let
$D=D(r)=(-r,r)^2$, and choose $r$ small enough so that $\bP^0(
Y_{t_0}\in D(r))< 1/20$. Let
$$ \tau_D = \inf\{ t>0: Y_t \not\in  D \}
\qq  \hbox{and} \qq  F(\lam) =\{ Y_{\tau_D} \in V(\lam) \} . $$
 Note that if $Y_{\tau_D} \not\in V(\lam)$ and $\tau_{V(\lam)}>t_0$ then
we have $ \tau_D\ge \tau_{V(\lam)}>t_0$. So,
$$ \bP^0(F(\tfrac13)^c) \le
 \bP^0(\tau_{V(\frac13)} \le t_0) + \bP^0( Y_{t_0} \in D) \le 1/10. $$
Note also that the events $F(\lam)$ are increasing in $\lam$.
}
\end{remark}

\begin{lemma} \label{Yharm}
There exists  sequence $x_n \to 0$ and
$\delta_n \to 0$ such that
 $\bigcup_{n=1}^\infty B(x_n, \delta_n) \subset V(\half)$ and
\[    \bP^y(F(\half)) \ge 7/10 \qq \hbox{ for } y \in \bigcup_{n=1}^\infty
 B(x_n, \delta_n). \]
\end{lemma}

\proof Let $h_Y(x)=\bP^x(F(\half))$. Then $h_Y(0)>9/10$
and $h_Y(Y_{t\wedge \tau_D})$ is
a right-continuous martingale. Using the right-continuity of $h_Y(Y)$
and $Y$, and the fact that $Y_t \in V(1/4)$ for all small times $t$,
we deduce that there exist $x_n \in V(1/4)$ with $x_n\to 0$ such that
$h_Y(x_n)\ge 8/10$ for all $n$.

Note that, since the coordinate processes $Y^i$ are symmetric
stable processes, each point $x \in \partial D$ is regular for $D^c$.
Therefore, if the event $F(\frac13)$ occurs for the process $Y(\omega)$
(with $Y_0=x_n$) then $F(\half)$ occurs, a.s., for the process
$u+Y$ for all sufficiently small $u$.
It follows that, for each $n$,
$$ \liminf_{y \to x_n} \bP^y(F(\half)) \ge 8/10. $$
We now take $\delta_n>0$ small enough so that
$B(x_n,\delta_n) \in V(\half)$ and $\bP^y(F(\half))\ge 7/10$
on $B(x_n,\delta_n)$.  \qed

\smallskip

We now define a second jump kernel $J_1$ follows.
Write $B=B(0,1)$.
If $x,y\in V$ we set
\[  J_1(x,y)= m\big(|x_1-y_1|,|x_2-y_2|\big)\mathbbm{1}_{(x-y \in B)} . \]
If $x,y\in V^c$ we set
\[  J_1(x,y)= m\big(|x_2-y_2|,|x_1-y_1|\big)\mathbbm{1}_{(x-y \in B)} . \]
If $x\in V$ and $y\in V^c$ or vice versa and $|x-y|\leq 1$, we define
\[  J_1(x,y)=\left( |x_1-y_1|^{-2-a} \wedge |x_2-y_2|^{-2-a} \right) \mathbbm{1}_{(x-y \in B)}. \]

It is easy to see that $J_1$ satisfies (A1)--(A3).

\medskip
\noindent{\bf Proof of Theorem \ref{Theo-counterexample}:}
Let $X=\{X_t, t\geq 0\}$ be the symmetric jump process associated
with the Dirichlet form
given by (\ref{Dirform})--(\ref{eqn:domain}) with $J_1$ as above
in place of $J$.
Note that if $x \in V$, then $J_1(x,y)=J_0(x,y)$ for $y \in V$, while
$J_1(x,y) \leq J_0(x,y)$ for $y\in V^c$. Thus $X$ makes as many
jumps within  $V$ as $Y$ does, but makes fewer jumps from $V$ to
$V^c$. This can be made more precise by using Remark \ref{meyer} to
construct $Y$ from $X$. Although when $X$ enters $V^c$, the positive
continuous additive functional $C=\{C_t, t\geq 0\}$ defined in
Remark \ref{meyer} becomes infinite immediately, we will only be
looking at time intervals $[0,\tau_V(X)]$, so this will not be an
issue for us.
In particular, all that we need is that $X_s=Y_s$
for $0\le s < \tau_V(Y)$.

Let $\sigma=\inf\{ t>0: X_t \not\in  D \}$, and set
$$ h(x)= \bP^x( X_{\sigma}\in V(\half)). $$

Then as $X=Y$ on $[0,\tau_V(Y))$, by a similar argument as that for
Lemma \ref{Yharm}, there exists $E:=\cup_{n\geq 1} B(x_n,
\delta_n)\subset V(\half )$, where $x_n\to 0$ and $\delta_n\to 0$,
such that
 $$h(x)=   \bP^x( X_{\sigma} \in V(\half ) )  \ge 7/10
 \qquad \hbox{for } x\in E.
 $$

Let $\Theta: \bR^2\to \bR^2$ be defined by
\begin{equation}\label{7Th}
\Theta((x_1,x_2))=(x_2,x_1).
\end{equation}
Then the reflection symmetry of the law of $X$ gives
$$ \bP^x( X_\sigma \in \Theta( V(\half)) )\ge 7/10, \qq
 x \in \Theta(E). $$
Hence we have
$$ h(x) \le 3/10, \qq x \in  \Theta(E). $$
Since $E$ and  $\Theta(E)$ are open sets and $0$
is an accumulation point of both sets, we deduce that $h$ is
not continuous at $0$.
\qed

\begin{corollary}\label{C7.1}
There exist a bounded continuous function $H$  and $t_0>0$ such that
$x\mapsto \E^x  \left[ H(X_{t_0}) \right]$ is not continuous at 0.
\end{corollary}

\proof Recall that  $Y$ is the L\'evy process with jump intensity
kernel $J_0$.
 By Remark \ref{R7.1} there exist $t_0>0$ and $r$ such
that
$$\P^0\left(Y_{t_0}\in V(1/3)\setminus D(r) \right)\geq \frac{9}{10}.
$$
Let $H$ be a continuous function bounded by $-1$ and 1 such that $H$
is 1 on $V(\frac13)\setminus D(r)$, $H\geq 0$ on $V(1)$, and
$H(x_2,x_1)=-H(x_1,x_2)$. Since $H$ is bounded below by $-1$, then
$\E^0 \left[ H(Y_{t_0})\right]\geq
\frac9{10}-\frac{1}{10}=\frac{8}{10}$. Since $Y$ is a L\'evy
process, it has the Feller property, and so $x\mapsto \E^x \left[
H(Y_{t_0})\right]$ is continuous. Therefore $\lim_{x\to 0} \E^x
\left[ H(Y_{t_0})\right]\geq \frac{8}{10}$.

As in the proof of Theorem \ref{Theo-counterexample}, there exist
sequences $x_n\to 0$, $\delta_n\to 0$ such that $E:=\cup_n B(x_n,
\delta_n)\subset V(1)$ and $\E^x\left[ H(X_{t_0})\right]\geq
\frac{7}{10}$ for $x\in E$. By the antisymmetry of $H$, $\E^x \left[
H(X_{t_0}) \right] \leq -\frac7{10}$ for $x\in \Theta (E)$, where
$\Theta$ is defined by (\ref{7Th}). We conclude $x\mapsto \E^x
\left[ H(X_{t_0}) \right]$ is not continuous at 0. \qed

\medskip
\begin{remark}\label{rem2} {\rm
\begin{description}

\item{(i)}  If we wish we can replace the double cone $V(\lam)$
by the half-space $H_-=\{ x_1\le 0\}$.
Use polar coordinates in $\bR^2$ and let $\theta_0 \in (0,\pi/2)$
satisfy $\tan \theta_0 = 1/\lam$. Let
$\psi: [0,2\pi] \to [0,2\pi]$ be strictly increasing and piecewise
linear with $\psi(0)=0$, $\psi(\theta_0)=\pi/2$,
$\psi(2\pi-\theta_0)=3\pi/2$ and $\psi(2\pi)=2\pi$.
Let $\Psi(r,\theta)=(r,\psi(\theta))$; note that
$\Psi(V(\lam)) \subset H_-$.

Let $Z=\Psi(X)$. Then $Z$ is a
symmetric jump process associated with a Dirichlet form $\sE'$ and
jump measure $n'$ which satisfies (A1)--(A3). By the construction of
$Z$ and using Remark \ref{R7.1} we see that there exist points $x_n \to 0$ with
$x_n \in H_-$, and $t_0>0$ such that,  for all $n$,
\begin{equation}\label{zhineq}
\bP^{x_n} \left( \tau^Z_{H_-} \geq  t_0 \right) \geq 9/10.
\end{equation}
If we time change $Z$ so that the associated Dirichlet form
is on $L^2(\bR^2,dx)$, the new process $\wt Z$ still
satisfies (\ref{zhineq}) with a different value of $t_0$.

 In particular we see that for the
processes considered in this paper, if $H$ is a half space, points
on $\partial H$ need not be regular for $H$.

\item{(ii)} In this example we needed $\al<\beta<1$ because points are
not polar for the symmetric stable process with $\alpha>1$. If we
look at similar constructions in higher dimensions such as $\bR^4=
\bR^2 \times \bR^2$ then it seems likely that we could construct a
similar example for any $0<\alpha<\beta <2$.

\item{(iii)} Note that (\ref{eq:m-choice}) allows for the choice of
$b=a+\eps$ for any value of $\eps >0$. Define
\begin{equation}
 \widetilde{m}(z_1,z_2) =
\begin{cases}
    \min\left\{|z_1|^{-a-2}, \log(\frac{3}{|z_2|})|z_2|^{-a-2} \right\}
  &\hbox{ if } |z_1|\vee |z_2| \le 1, \\
    0                         &\hbox{ if } |z_1|\vee |z_2| > 1. \\
\end{cases}
\end{equation}
Then $\widetilde{m}$ is very similar to $m$ in (\ref{eq:m-choice}).
Now, proceed as in the example above but this time construct $J_1$
with the help of  $\widetilde{m}$ instead of $m$. Then it is shown
as a corollary in \cite{HuKa05} that the Dirichlet form corresponds
to a Feller semigroup and that harmonic functions are continuous
satisfying certain {\em a priori} estimates. Therefore when trying
to construct discontinuous harmonic functions, one cannot modify our
class of examples much.
   \end{description} }
\end{remark}

\begin{proposition}\label{MGproblem}
With $J$ as above, the martingale
problem for the operator
\[
\sL f(x)= \int_{\R^2} ( f(x+h)-f(x) ) J_1(x,x+h)\, dh
\]
acting on $C^2_c$ functions
is not well-posed.
\end{proposition}

\proof The function $\sL f$ is bounded by the $C^2$-norm of $f$ when
$f\in C^2_c(\R^2)$ by (A1)--(A3). Section 3 of  \cite{Bas88}
(trivially modified to handle the case of dimensions larger than
one) shows that $\{\P^x: x\in \R^2\}$ is tight.

 We claim that any
subsequential weak limit point of $\P^x$ as $x\to 0$ is a solution
to the martingale problem for $\sL$ started at 0.  If $\P$ is such a
limit point, it is easy to see that $\P(X_0=0)=1$. If $f\in C^2_c$
and $r_1\leq \cdots \leq r_n\leq s\leq t$ and $g_1, \ldots, g_n$ are
continuous functions with compact support, then
\begin{align}
\E^{x_n}\Big[ \Big\{ &f(X_t)-f(X_0)-\int_0^t \sL f(X_u)\, du\Big\} Y\Big]
\label{MGtest}\\
&= \E^{x_n}\Big[ \Big\{ f(X_s)-f(X_0)-\int_0^s \sL f(X_u)\, du\Big\} Y\Big],\nn
\end{align}
where $Y=\prod_{i=1}^n g_i(X_{r_i})$. Since $X$ has no jumps larger
than 1, we see that if $f\in C^2_c (\R^2) $, then $\sL f$ will be
zero at points that are a distance more than one from the support of
$f$. Therefore $\sL f$ also has compact support. Let $\eps>0$ and
let $I$ be a continuous function that equals $\sL f$ except on a set
$A$ of Lebesgue measure at most $\eps$. Since $f(X_t)Y$, $f(Y_s)Y$,
$f(X_0)Y$, $Y\int_0^t I(X_u)\, du$, and $Y\int_0^s I(X_u)\, du$ are
continuous functionals of the path, their expectation under
$\P^{x_n}$ converges to the corresponding expectation under $\P$. We
have the estimate
$$\E^{X_n}\int_0^t |\sL f(X_u)-I(X_u)|\, du
\leq c_1 t_1 +c_1\E^{x_n}\left[\int_{t_1}^t \mathbbm{1}_A(X_u)\, du\right]
\leq c_1t_1+ c_2t_1^{-d/\al}|A|$$
using Theorem \ref{rough-above-lem}. If we set $t_1=|A|^{\al/(\al+d)}$,
we obtain the upper bound $c_3|A|^{\al/(\al+d)}$.
A limit argument yields  the same bound for the  expectation with
respect to $\P$.
Therefore
\begin{align*}
\limsup_{n\to \infty} \Bigl|&\E\Big[
\Big\{ f(X_t)-f(X_0)-\int_0^t \sL f(X_u)\, du\Big\} Y\Big]\\
&~~~~~~~~~~~~~~~
- \E\Big[ \Big\{ f(X_s)-f(X_0)-\int_0^s \sL f(X_u)\, du\Big\} Y\Big]\Bigr|\\
&\leq c_4 \eps^{\al/(\al+d)}.
\end{align*}
Since $\eps$ is arbitrary, we have
(\ref{MGtest}) with $\E^{x_n}$ replaced by $\E$. This shows that
$\P$ is a solution to the martingale problem started at 0.

Suppose now that we had uniqueness to the martingale problem for
$\sL$ started at 0. We conclude $\P^x\to \P^0$ as $x\to 0$. In
particular, if $H$ is a bounded continuous function, $\E^x
\left[H(X_t) \right]\to \E^0 \left[H(X_t) \right]$ as $x\to 0$ for
all $t>0$. But this contradicts Corollary \ref{C7.1}. \qed

\bigskip
\bigskip

\noindent{\bf Martin T. Barlow}\\
Department of Mathematics\\
University of  British Columbia \\
Vancouver, BC, V6T 1Z2,  Canada\\
{\texttt barlow@math.ubc.ca}\\

\medskip
\noindent{\bf Richard. F. Bass}\\
Department of Mathematics\\
University of  Connecticut \\
Storrs, CT 06269-3009, USA\\
{\texttt bass@math.uconn.edu} \\

\medskip
\noindent{\bf Zhen-Qing Chen}\\
Department of Mathematics\\
University of Washington\\
Seattle, WA 98195, USA\\
{\texttt zchen@math.washington.edu}\\

\medskip
\noindent{\bf Moritz Kassmann}\\
Institut f\"{u}r Angewandte Mathematik\\
Universit\"{a}t Bonn \\
Beringstrasse 6\\
D-53115 Bonn, Germany\\
{\texttt kassmann@iam.uni-bonn.de}

\end{document}